\documentclass[notitlepage,leqno,10pt]{article}

\usepackage{latexsym}
\usepackage{amsmath}
\usepackage{amsfonts}
\usepackage{amssymb}
\usepackage{graphicx}
\renewcommand{\a }{\alpha }

\renewcommand{\d}{\delta }
\newcommand{\D }{\Delta }

\newcommand{\e }{\varepsilon }
\newcommand{\g }{\gamma}
\renewcommand{\i }{\iota}

\renewcommand{\l }{\lambda }
\renewcommand{\L }{\Lambda }

\newcommand{\n }{\nabla }
\newcommand{\var }{\varphi }

\newcommand{\s }{\sigma }
\newcommand{\Sig }{\Sigma}

\renewcommand{\th }{\theta }

\newcommand{\ov}{\overline}
\newcommand{\intbar}{\mathop{\int\makebox(-13.5,0){\rule[4pt]{.7em}{0.3pt}}%
\kern-6pt}\nolimits}

\newcommand{\be}{\begin{equation}}
\newcommand{\ee}{\end{equation}}
\newenvironment{pf}{\noindent{\sc Proof}.\enspace}{\rule{2mm}{2mm}\medskip}
\newenvironment{pfn}{\noindent{\sc Proof}}{\rule{2mm}{2mm}\medskip}
\newenvironment{pfnb}{\noindent{\bf Proof}}{\rule{2mm}{2mm}\medskip}

\newcommand{\R}{\mathbb{R}}

\newcommand{\N}{\mathbb{N}}

\author{Zindine DJADLI$^{a}$ and Andrea MALCHIODI$^{b}$}

\date{}

\title{Existence of conformal metrics with constant $Q$-curvature}

\begin{document}

\newtheorem{lem}{Lemma}[section]
\newtheorem{pro}[lem]{Proposition}
\newtheorem{thm}[lem]{Theorem}
\newtheorem{rem}[lem]{Remark}
\newtheorem{cor}[lem]{Corollary}
\newtheorem{df}[lem]{Definition}

\maketitle

\begin{center}

{\small $^a$ Universit\'e de Cergy-Pontoise - D\'epartement de
Math\'ematiques

Site de Saint-Martin, 2 avenue Adolphe Chauvin

 F95302
Cergy-Pontoise Cedex - France

\

\noindent $^b$ SISSA, via Beirut 2-4, 34014 Trieste, Italy.}

\end{center}

\footnotetext[1]{E-mail addresses: zindine.djadli@math.u-cergy.fr,
malchiod@sissa.it}

\

\

\noindent {\sc abstract}. Given a compact four dimensional manifold,
we prove existence of conformal metrics with constant $Q$-curvature
under generic assumptions. The problem amounts to solving a
fourth-order nonlinear elliptic equation with variational structure.
Since the corresponding Euler functional is in general unbounded
from above and from below, we employ topological methods and min-max
schemes, jointly with the compactness result of \cite{mal}.
\begin{center}

\bigskip\bigskip


\bigskip

\centerline{\bf AMS subject classification: 35B33, 35J35, 53A30,
53C21}

\end{center}

\section{Introduction}\label{s:in}

In recent years, much attention has been devoted to the study of
partial differential equations on manifolds, in order to understand
some connections between analytic and geometric properties of these
objects.

A basic example is the Laplace-Beltrami operator on a compact
surface $(\Sig,g)$. Under the conformal change of metric
$\tilde{g} = e^{2w}g$, we have
\begin{equation}\label{eq:lapl}
    \Delta_{\tilde{g}} = e^{-2w} \Delta_g; \qquad \qquad -
    \Delta_g w + K_g = K_{\tilde{g}} e^{2w},
\end{equation}
where $\D_g$ and $K_g$ (resp. $\D_{\tilde{g}}$ and
$K_{\tilde{g}}$) are the Laplace-Beltrami operator and the Gauss
curvature of $(\Sig, g)$ (resp. of $(\Sig, \tilde{g})$). From the
above equation one recovers in particular the conformal invariance
of $\int_\Sig K_g dV_g$, which is related to the topology of
$\Sig$ through the Gauss-Bonnet formula
\begin{equation}
\int_\Sig K_g dV_g = 2\pi \chi(\Sig),
\end{equation}
where $\chi(\Sig)$ is the Euler characteristic of $\Sig$. Of
particular interest is the classical {\em Uniformization Theorem},
which asserts that every compact surface carries a (conformal)
metric with constant curvature.

\

\noindent On four dimensional manifolds there exists a conformally
covariant operator, the Paneitz operator, which enjoys analogous
properties to the Laplace-Beltrami operator on surfaces, and to
which is associated a natural concept of curvature. This operator,
introduced by Paneitz, \cite{p1}, \cite{p2}, and the corresponding
$Q$-curvature, introduced in \cite{br2}, are defined in terms of
Ricci tensor $Ric_g$ and scalar curvature $R_g$ of the manifold
$(M,g)$ as
\begin{equation}
P_g (\varphi) = \Delta_g^2 \varphi + div_g \left(\frac 23 R_g g -
2 Ric_g \right) d \varphi;
\end{equation}
\begin{equation}\label{eq:Q}
    Q_g = - \frac{1}{12} \left( \D_g R_g - R_g^2 + 3 |Ric_g|^2
    \right),
\end{equation}
where $\varphi$ is any smooth function on $M$. The behavior (and
the mutual relation) of $P_g$ and $Q_g$ under a conformal change
of metric $\tilde{g} = e^{2w}g$ is given by
\begin{equation}\label{eq:confP}
P_{\tilde{g}} = e^{- 4 w} P_g; \qquad  \qquad P_g w + 2 Q_g =
    2 Q_{\tilde{g}} e^{4 w}.
\end{equation}
Apart from the analogy with \eqref{eq:lapl}, we have an extension
of the Gauss-Bonnet formula which is the following
\begin{equation}\label{eq:GB4}
    \int_M \left( Q_g + \frac{|W_g|^2}{8} \right) dV_g = 4 \pi^2
    \chi(M),
\end{equation}
where $W_g$ denotes the Weyl tensor of $(M,g)$ and $\chi(M)$ the
Euler characteristic. In particular, since $\left\vert W_g
\right\vert^2dV_g$ is a pointwise conformal invariant, it follows
that the integral of $Q_g$ over $M$ is also a conformal invariant,
which is usually denoted with the symbol
\begin{equation}\label{eq:kp}
    k_P = \int_M Q_g dV_g.
\end{equation}
We refer for example to the survey \cite{cy99} for more details.

To mention some first geometric properties of $P_g$ and $Q_g$, we
discuss some results of Gursky, \cite{g2} (see also \cite{g}). If a
manifold of non-negative Yamabe class $Y(g)$ (this means that there
is a conformal metric with non-negative constant scalar curvature)
satisfies $k_P \geq 0$, then the kernel of $P_g$ are only the
constants and $P_g \geq 0$, namely $P_g$ is a non-negative operator.
If in addition $Y(g) > 0$, then the first Betti number of $M$
vanishes, unless $(M,g)$ is conformally equivalent to a quotient of
$S^3 \times \R$. On the other hand, if $Y(g) \geq 0$, then $k_P \leq
8 \pi^2$, with the equality holding if and only if $(M,g)$ is
conformally equivalent to the standard sphere.

\

\noindent As for the Uniformization Theorem, one can ask whether
every four-manifold $(M,g)$ carries a conformal metric $\tilde{g}$
for which the corresponding $Q$-curvature $Q_{\tilde{g}}$ is a
constant. Writing $\tilde{g} = e^{2 w} g$, by \eqref{eq:confP} the
problem amounts to finding a solution of the equation
\begin{equation}\label{eq:Qc}
    P_g w + 2 Q_g = 2 \ov{Q} e^{4 w},
\end{equation}
where $\ov{Q}$ is a real constant. By the regularity results in
\cite{uv}, critical points of the following functional
\begin{equation}\label{eq:II}
    II(u) = \langle P_g u, u \rangle + 4 \int_M Q_g u dV_g - k_P
    \log \int_M e^{4u} dV_g; \qquad u \in H^{2}(M),
\end{equation}
which are weak solutions of \eqref{eq:Qc}, are also strong
solutions. Here $H^2(M)$ is the space of functions on $M$ which are
of class $L^2$, together with their first and second derivatives,
and the symbol $\langle P_g u, v \rangle$ stands for
\begin{equation}\label{eq:Pguv}
    \langle P_g u, v \rangle = \int_M \left( \D_g u \D_g v + \frac 23
    R_g \n_g u \cdot \n_g v - 2 (Ric_g \n_g u, \n_g v) \right) dV_g
    \qquad \hbox{ for } u, v \in H^2(M).
\end{equation}

Problem \eqref{eq:Qc} has been solved in \cite{cy95} for the case in
which $ker P_g = \R$, $P_g$ is a non-negative operator and $k_P < 8
\pi^2$. By the above-mentioned result of Gursky, sufficient
conditions for these assumptions to hold are that $Y(g) \geq 0$ and
that $k_P \geq 0$ (and $(M,g)$ is not conformal to the standard
sphere). Notice that if $Y(g) \geq 0$ and $k_P = 8 \pi^2$, then
$(M,g)$ is conformally equivalent to the standard sphere and clearly
in this situation \eqref{eq:Qc} admits a solution. More general
conditions for the above hypotheses to hold have been obtained by
Gursky and Viaclovsky in \cite{gv}. Under the assumptions in
\cite{cy95}, by the Adams inequality
$$
  \log \int_M e^{4(u - \ov{u})} dV_g \leq \frac{1}{8 \pi^2}
    \langle P_g u, u \rangle + C, \qquad \quad u \in H^2(M),
$$
where $\ov{u}$ is the average of $u$ and where $C$ depends only on
$M$, the functional $II$ is bounded from below and coercive, hence
solutions can be found as global minima. The result in \cite{cy95}
has also been extended in \cite{b1} to higher-dimensional manifolds
(regarding higher-order operators and curvatures) using a geometric
flow.

The solvability of \eqref{eq:Qc}, under the above hypotheses, has
been useful in the study of some conformally invariant fully
non-linear equations, as is shown in \cite{cgy}. Some remarkable
geometric consequences of this study, given in \cite{cgyann},
\cite{cgy}, are the following. If a manifold of positive Yamabe
class satisfies $\int_M Q_g dV_g
> 0$, then there exists a conformal metric with positive Ricci
tensor, and hence $M$ has finite fundamental group. Furthermore,
under the additional quantitative assumption $\int_M Q_g dV_g >
\frac 18 \int_M |W_g|^2 dV_g$, $M$ must be diffeomorphic to the
standard four-sphere or to the standard projective space. Finally,
we also point out that the Paneitz operator and the $Q$-curvature
(together with their higher-dimensional analogues, see \cite{br},
\cite{br2}, \cite{fg}, \cite{gjms}) appear in the study of
Moser-Trudinger type inequalities, log-determinant formulas and the
compactification of locally conformally flat manifolds, \cite{bo},
\cite{bcy}, \cite{cqy}, \cite{cqy2}, \cite{cy95}.

\

\noindent We are interested here in extending the {\em
uniformization} result in \cite{cy95}, namely to find solutions of
\eqref{eq:Qc} under more general assumptions. Our result is the
following.

\begin{thm}\label{th:ex}
Suppose $ker \; P_g = \{constants\}$, and assume that $k_P \neq 8 k
\pi^2$ for $k = 1, 2, \dots$. Then $(M,g)$ admits a conformal metric
with constant $Q$-curvature.
\end{thm}

\begin{rem} (a) Our assumptions are conformally invariant and generic,
so the result applies to a large class of four manifolds, and in
particular to some manifolds of negative curvature or negative
Yamabe class. Note that, in view of \cite{g2}, it is not clear
whether or not a manifold of negative Yamabe class satisfies the
assumptions of the result in \cite{cy95}. For example, products of
two negatively-curved surfaces might have total $Q$-curvature
greater than $8 \pi^2$, see \cite{dm1}.

(b) Under these assumptions, imposing the volume normalization
$\int_M e^{4u} dV_g = 1$, the set of solutions (which is non-empty)
is bounded in $C^m(M)$ for any integer $m$, by Theorem 1.3 in
\cite{mal}, see also \cite{dr}.

(c) Theorem \ref{th:ex} does NOT cover the case of locally
conformally flat manifolds with positive and even Euler
characteristic, by \eqref{eq:GB4}.
\end{rem}

\

\noindent Our assumptions include those made in \cite{cy95} and
one (or both) of the following two possibilities
\begin{equation}\label{eq:kp3}
    k_P \in (8 k \pi^2, 8 (k+1) \pi^2), \quad \hbox{ for some } k \in \N;
\end{equation}
\begin{equation}\label{eq:kp4}
    P_g \hbox{ possesses $\ov{k}$ (counted with
    multiplicity)
    negative eigenvalues}.
\end{equation}
In these cases the functional $II$ is unbounded from below, and
hence it is necessary to find extremals which are possibly saddle
points. This is done using a new min-max scheme, which we are going
to describe below, depending on $k_P$ and the spectrum of $P_g$ (in
particular on the number of negative eigenvalues $\ov{k}$, counted
with multiplicity). By classical arguments, the scheme yields a {\em
Palais-Smale sequence}, namely a sequence $(u_l)_l \subseteq H^2(M)$
satisfying the following properties
\begin{equation}\label{eq:PS}
    II(u_l) \to c \in \R; \qquad II'(u_l) \to 0 \qquad \qquad
    \hbox{ as } l \to + \infty.
\end{equation}
We can also assume that such a sequence $(u_l)_l$ satisfies the
volume normalization
\begin{equation}\label{eq:noul}
    \int_M e^{4 u_l} dV_g = 1 \qquad \qquad \hbox{ for all } l.
\end{equation}
This is always possible since the functional $II$ is invariant under
the transformation $u \mapsto u + a$, where $a$ is any real
constant. Then, to achieve existence, one should prove for example
that $(u_l)_l$ is bounded, or to prove a similar compactness
criterion.

In order to do this, we apply a procedure from \cite{str}, used in
\cite{djlw}, \cite{jt}, \cite{st}. For $\rho$ in a neighborhood of
$1$, we define the functional $II_{\rho} : H^2(M) \to \R$ by
$$
  II_\rho(u) = \langle P_g u, u \rangle + 4 \rho \int_M Q_g d V_g
  - 4 \rho k_P \log \int_M e^{4 u} d V_g, \qquad u \in H^2(M),
$$
whose critical points give rise to solutions of the equation
\begin{equation}\label{eq:mod}
    P_g u + 2 \rho Q_g = 2 \rho k_P e^{4 u} \qquad \hbox{ in } M.
\end{equation}
One can then define the min-max scheme for different values of
$\rho$ and prove boundedness of some Palais-Smale sequence for
$\rho$ belonging to a set $\L$ which is dense in some neighborhood
of $1$, see Section \ref{s:proof}. This implies solvability of
\eqref{eq:mod} for $\rho \in \L$. We then apply the following result
from \cite{mal}, with $Q_l = \rho_l Q_g$, where $(\rho_l)_l
\subseteq \L$ and $\rho_l \to 1$.

\begin{thm}\label{th:bd} (\cite{mal})
Suppose $ker \; P_g = \{constants\}$ and that $(u_l)_l$ is a
sequence of solutions of
\begin{equation}\label{eq:pl}
    P_g u_l + 2 Q_l = 2 k_l e^{4 u_l} \qquad \hbox{ in } M,
\end{equation}
satisfying \eqref{eq:noul}, where $k_l = \int_M Q_l d V_g$, and
where $Q_l \to Q_0$ in $C^0(M)$. Assume also that
\begin{equation}\label{eq:kp2}
    k_0 := \int_M Q_0 dV_g \neq 8 k \pi^2 \qquad \qquad
    \hbox{ for } k = 1, 2, \dots.
\end{equation}
Then $(u_l)_l$ is bounded in $C^\a(M)$ for any $\a \in (0,1)$.
\end{thm}

\

\noindent We are going to give now a brief description of the scheme
and an heuristic idea of its construction. We describe it for the
functional $II$ only, but the same considerations hold for $II_\rho$
if $|\rho -1|$ is sufficiently small. It is a standard method in
critical point theory to find extrema by looking at the difference
of topology between sub or superlevels of the functionals. In our
specific case we investigate the structure of the sublevels $\{ II
\leq - L\}$, where $L$ is a large positive number. Looking at the
form of the functional $II$, see \eqref{eq:II}, one can image two
ways for attaining large negative values.

The first, assuming \eqref{eq:kp3}, is by bubbling. In fact, for a
given point $x \in M$ and for $\l > 0$, consider the following
function
$$
  \var_{\l,x}(y) = \log \left( \frac{2 \l}{1 + \l^2 dist(y,x)^2}
  \right),
$$
where $dist(\cdot,\cdot)$ denotes the metric distance on $M$. Then
for $\l$ large one has $e^{4 \var_{\l,x}} \simeq \d_x$ (the Dirac
mass at $x$), where $e^{4 \var_{\l,x}}$ represents the volume
density of a four sphere attached to $M$ at the point $x$, and one
can show that $II(\var_{\l,x}) \to - \infty$ as $\l \to + \infty$.
Similarly, for $k$ given in \eqref{eq:kp3} and for $x_1, \dots, x_k
\in M$, $t_1, \dots, t_k \geq 0$, it is possible to construct an
appropriate function $\var$ of the above form (near each $x_i$) with
$e^{4 \var} \simeq \sum_{i = 1}^k t_i \d_{x_i}$, and on which $II$
still attains large negative values. Precise estimates are given in
Section \ref{s:test} and in the appendix. Since $II$ stays invariant
if $e^{4 \var}$ is multiplied by a constant, we can assume that
$\sum_{i=1}^k t_i = 1$. On the other hand, if $e^{4 \var}$ is
concentrated at $k+1$ distinct points of $M$, it is possible to
prove, using an improved Moser-Trudinger inequality from Section
\ref{s:pr}, that $II(\var)$ cannot attain large negative values
anymore, see Lemmas \ref{l:imprc} and \ref{l:II<-M}. From this
argument we see that one is led naturally to consider the family
$M_k$ of elements $\sum_{i=1}^k t_i \d_{x_i}$ with $(x_i)_i
\subseteq M$, and $\sum_{i=1}^k t_i = 1$, known in literature as the
{\em formal set of barycenters of} $M$ {\em of order} $k$, which we
are going to discuss in more detail below.

The second way to attain large negative values, assuming
\eqref{eq:kp4}, is by considering the negative-definite part of the
quadratic form $\langle P_g u, u \rangle$. Letting $V \subseteq
H^2(M)$ denote the direct sum of the eigenspaces of $P_g$
corresponding to negative eigenvalues, the functional $II$ will tend
to $- \infty$ on the boundaries of large balls in $V$, namely
boundaries sets homeomorphic to the unit ball $B^{\ov{k}}_1$ of
$\R^{\ov{k}}$.

\

Having these considerations in mind, we will use for the min-max
scheme a set, denoted by $A_{k,\ov{k}}$, which is constructed using
some contraction of the product $M_k \times B^{\ov{k}}_1$, see
formula \eqref{eq:akovk} and the figure in Section \ref{s:pr} (when
$k_P < 8 \pi^2$, we just take the sphere $S^{\ov{k}-1}$ instead of
$A_{k,\ov{k}}$). It is possible indeed to map (non-trivially) this
set into $H^2(M)$ in such a way that the functional $II$ on the
image is close to $- \infty$, see Section \ref{s:test}. On the other
hand, it is also possible to do the opposite, namely to map
appropriate sublevels of $II$ into $A_{k,\ov{k}}$, see Section
\ref{s:msub}. The composition of these two maps turns out to be
homotopic to the identity on $A_{k,\ov{k}}$ (which is
non-contractible by Corollary \ref{c:akk}) and therefore they are
both topologically non-trivial.

\

\noindent Some comments are in order. For the case $k = 1$ and
$\ov{k} = 0$, which is presented in \cite{dm1}, the min-max scheme
is similar to that used in \cite{djlw}, where the authors study a
mean field equation depending on a real parameter $\l$ (and prove
existence for $\l \in (8 \pi, 16 \pi)$). Solutions for large values
of $\l$ have been obtained recently by Chen and Lin, \cite{cl1},
\cite{cl2}, using blow-up analysis and degree theory. See also the
papers \cite{li}, \cite{ln}, \cite{st} and references therein for
related results. The construction presented in this paper has been
recently used in Djadli \cite{dj} to study this problem as well when
$\l \neq 8k\pi$ and without any assumption on the topology on the
surface. Our method has also been employed by Malchiodi and Ndiaye
\cite{mn} for the study of the $2 \times 2$ Toda system.

The set of barycenters $M_k$ (see Subsection \ref{ss:prop} for more
comments or references) has been used crucially in literature for
the study of problems with lack of compactness, see \cite{bah},
\cite{bc}. In particular, for Yamabe-type equations (including the
Yamabe equation and several other applications), it has been used to
understand the structure of the {\em critical points at infinity}
(or asymptotes) of the Euler functional, namely the way compactness
is lost through a pseudo-gradient flow. Our use of the set $M_k$,
although the map $\Phi$ of Section \ref{s:test} presents some
analogies with the Yamabe case, is of different type since it is
employed to reach low energy levels and not to study critical points
at infinity. As mentioned above, we consider a projection onto the
$k$-barycenters $M_k$, but starting only from functions in $\{ II
\leq - L \}$, whose concentration behavior is not as clear as that
of the asymptotes for the Yamabe equation. Here also a technical
difficulty arises. The main point is that, while in the Yamabe case
all the coefficients $t_i$ are bounded away from zero, in our case
they can be arbitrarily small, and hence it is not so clear what the
choice of the points $x_i$ and the numbers $t_i$ should be when
projecting. Indeed, when $k > 1$ $M_k$ is not a smooth manifold but
a {\em stratified set}, namely union of sets of different dimensions
(the maximal one is $5k - 1$, when all the $x_i$'s are distinct and
all the $t_i$'s are positive). To construct a continuous global
projection takes us some work, and this is done in Section
\ref{s:msub}.

\

\noindent The cases which are not included in Theorem \ref{th:ex}
should be more delicate, especially when $k_P$ is an integer
multiple of $8 \pi^2$. In this case non-compactness is expected, and
the problem should require an asymptotic analysis as in \cite{bah},
or a fine blow-up analysis as in \cite{li}, \cite{cl1}, \cite{cl2}.
Some blow-up behavior on open flat domains of $\R^4$ is studied in
\cite{ars}.

A related question in this context arises for the standard sphere
($k_P = 8 \pi^2$), where one could ask for the analogue of the {\em
Nirenberg's problem}. Precisely, since the $Q$-curvature of the
standard metric is constant, a natural problem is to deform the
metric conformally in such a way that the curvature becomes a given
function $f$ on $S^4$. Equation \eqref{eq:Qc} on the sphere admits a
non-compact family of solutions (classified in \cite{cy97}), which
all arise from conformal factors of M\"obius transformations. In
order to tackle this loss of compactness, usually finite-dimensional
reductions of the problem are used, jointly with blow-up analysis
and Morse theory. Some results in this direction are given in
\cite{b2}, \cite{ms} and \cite{wx} (see also references therein for
results on the Nirenberg's problem on $S^2$).

\

\noindent The structure of the paper is the following. In Section
\ref{s:pr} we collect some notation and preliminary results, based
on an improved Moser-Trudinger type inequality. We also introduce
the set $A_{k,\ov{k}}$ used to perform the min-max construction. In
Section \ref{s:msub} then we show how to map the sublevels $\{ II
\leq - L \}$ into $A_{k,\ov{k}}$. We begin by analyzing some
properties of the $k$-barycenters as a stratified set, in order to
understand the component of the projection involving the set $M_k$,
which is the most delicate. Then we turn to the construction of the
global map. In Section \ref{s:test} we show how to embed
$A_{k,\ov{k}}$ into the sublevel $\{ II \leq - L \}$ for $L$
arbitrarily large. This requires long and delicate estimates, some
of which are carried out in the appendix (which also contains other
technical proofs). Finally in Section \ref{s:proof} we prove Theorem
\ref{th:ex}, defining a min-max scheme based on the construction of
$A_{k,\ov{k}}$, solving the modified problem \eqref{eq:mod}, and
applying Theorem \ref{th:bd}.

\

\noindent An announcement of the present results is given in the
preliminary note \cite{dm1}.

\

\begin{center}

{\bf Acknowledgements}

\end{center}

\noindent We thank A.Bahri for explaining us the proof of Lemma
\ref{l:nonc}. This work was started when the authors were visiting
IAS in Princeton, and continued during their stay at IMS in
Singapore. A.M. worked on this project also when he was visiting ETH
in Z\"urich and Laboratoire Jacques-Louis Lions in Paris. They are
very grateful to all these institutions for their kind hospitality.
Z.D. has been supported by the {\em ACI-Jeunes Chercheurs :
M´etriques privil\'egi\'ees sur les vari\'et\'es \`a bord
2003/2006}. A.M. has been supported by M.U.R.S.T. under the national
project {\em Variational methods and nonlinear differential
equations}, and by the European Grant ERB FMRX CT98 0201.

\section{Notation and preliminaries}\label{s:pr}

In this section we fix our notation and we recall some useful known
facts. We state in particular an inequality of Moser-Trudinger type
for functions in $H^2(M)$, an improved version of it and some of its
consequences.

The symbol $B_r(p)$ denotes the metric ball of radius $r$ and center
$p$, while $dist(x,y)$ stands for the distance between two points
$x, y \in M$. $H^2(M)$ is the Sobolev space of the functions on $M$
which are in $L^2(M)$ together with their first and second
derivatives. The symbol $\| \cdot \|$ will denote the norm of
$H^2(M)$. If $u \in H^2(M)$, $\ov{u} = \frac{1}{|M|} \int_M u dV_g$
stands for the average of $u$. For $l$ points $x_1, \dots, x_l \in
M$ which all lie in a small metric ball, and for $l$ non-negative
numbers $\a_1, \dots, \a_l$, we will consider {\em convex
combinations} of the form $\sum_{i=1}^l \a_i x_i$, $\a_i \geq 0$,
$\sum_i \a_i = 1$. To do this, we can consider the embedding of $M$
into some $\R^n$ given by Whitney's theorem, take the convex
combination of the images of the points $(x_i)_i$, and project it
onto the image of $M$ (which we identify with $M$ itself). If
$dist(x_i,x_j) < \xi$ for $\xi$ sufficiently small, $i, j = 1,
\dots, l$, then this operation is well-defined and moreover we have
$dist\left(x_j, \sum_{i=1}^l \a_i x_i \right) < 2 \xi$ for every $j
= 1, \dots, l$. Note that these elements are just points, not to be
confused with the formal barycenters $\sum t_i \d_{x_i}$.

Large positive constants are always denoted by $C$, and the value of
$C$ is allowed to vary from formula to formula and also within the
same line. When we want to stress the dependence of the constants on
some parameter (or parameters), we add subscripts to $C$, as $C_\d$,
etc.. Also constants with subscripts are allowed to vary.

Since we allow $P_g$ to have negative eigenvalues, we denote by $V
\subseteq H^2(M)$ the direct sum of the eigenspaces corresponding to
negative eigenvalues of $P_g$. The dimension of $V$, which is
finite, is denoted by $\ov{k}$, and since $ker P_g = \R$, we can
find a basis of eigenfunctions $\hat{v}_1, \dots, \hat{v}_{\ov{k}}$
of $V$ (orthonormal in $L^2(M)$) with the properties
\begin{equation}\label{eq:hatv1k}
    P_g \hat{v}_i = \l_i \hat{v}_i, \quad i = 1, \dots, \ov{k};
  \qquad \int_M \hat{v}_i^2 dV_g = 1; \qquad \l_1 \leq \l_2 \leq
  \dots \leq \l_{\ov{k}} < 0 < \l_{\ov{k}+1} \leq \dots,
\end{equation}
where the $\l_i$'s are the eigenvalues of $P_g$ counted with
multiplicity. From \eqref{eq:hatv1k}, since $P_g$ has a divergence
structure, it follows immediately that $\int_M \hat{v}_i dV_g = 0$
for $i = 1, \dots, \ov{k}$. We also introduce the following
positive-definite (on the space of functions orthogonal to the
constants) pseudo-differential operator $P^+_g$
\begin{equation}\label{eq:Pg+}
  P^+_g u = P_g u - 2 \sum_{i=1}^{\ov{k}} \l_i \left( \int_M u
  \hat{v}_i dV_g \right) \hat{v}_i.
\end{equation}
Basically, we are reversing the sign of the negative eigenvalues of
$P_g$.

\

\noindent Now we define the set $A_{k,\ov{k}}$ to be used in the
existence argument, where $k$ is as in \eqref{eq:kp3}, and where
$\ov{k}$ is as in \eqref{eq:hatv1k}. We let $M_k$ denote the family
of formal sums
\begin{equation}\label{eq:Mk}
    M_k = \sum_{i=1}^k t_i \d_{x_i}; \qquad \qquad t_i \geq 0,
    \sum_{i=1}^k t_i = 1; \quad x_i \in M,
\end{equation}
endowed with the weak topology of distributions. This is known in
literature as the {\em formal set of barycenters} of $M$ (of order
$k$), see \cite{bah}, \cite{bc}, \cite{bred}. We stress that this
set is NOT the family of convex combinations of points in $M$ which
is introduced at the beginning of the section. To carry out some
explicit computations, we will use on $M_k$ the metric given by
$C^1(M)^*$, which induces the same topology, and which will be
denoted by $dist(\cdot, \cdot)$.

\begin{figure}[h]
\begin{center}
\includegraphics[width=13cm]{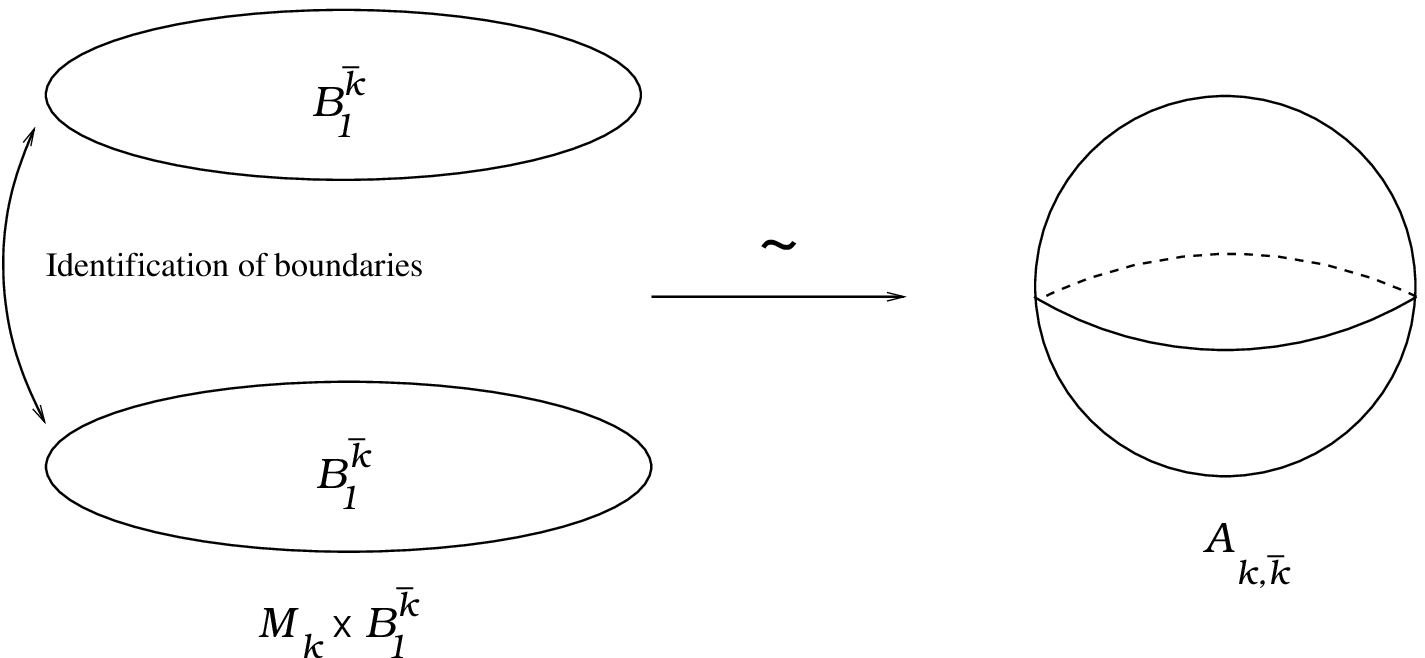}
\caption{the equivalence relation \; $\widetilde{}$} \label{fig1}
\end{center}
\end{figure}

Then, recalling that $\ov{k}$ is the number of negative eigenvalues
of $P_g$, we consider the unit ball $B_{1}^{\ov{k}}$ in
$\R^{\ov{k}}$, and we define the set
\begin{equation}\label{eq:akovk}
    A_{k,\ov{k}} = \widetilde{M_k \times B_{1}^{\ov{k}}},
\end{equation}
where the notation $\; \widetilde{} \;$ means that $M_k \times
\partial B_{1}^{\ov{k}}$ is identified with $\partial
B_{1}^{\ov{k}}$, namely $M_k \times \{y\}$, for every fixed $y \in
\partial B_{1}^{\ov{k}}$, is collapsed to a single point. In Figure
\ref{fig1} below we illustrate this collapsing drawing, for
simplicity, $M_k$ as a couple of points. When $k_P < 8 \pi^2$ and
$\ov{k} \geq 1$, we will perform the min-max argument just by using
the sphere $S^{\ov{k}-1}$.

\subsection{Some improved Adams inequalities}

\noindent In this subsection we give some improvements of the Adams
inequality (see \cite{ada} and \cite{cy95}) and in particular we
consider the possibility of dealing with operators $P_g$ possessing
negative eigenvalues. The following Lemma is proved in \cite{mal}
using a modification of the arguments in \cite{cy95}, which in turn
extend to the Paneitz operator some previous embeddings due to Adams
involving the operator $\D^m$ in flat domains.

\begin{lem}\label{l:adaneg} (\cite{mal})
Suppose $ker P_g = \{constants\}$, let $V$ be the direct sum of the
eigenspaces corresponding to negative eigenvalues of $P_g$, and let
$P_g^+$ be defined in \eqref{eq:Pg+}. Then there exists a constant
$C$ such that for all $u \in H^2(M)$
\begin{equation}\label{eq:adaneg0}
    \int_M e^{\frac{32 \pi^2 (u - \ov{u})^2}{\langle P_g^+ u, u \rangle}}
    dV_g \leq C.
\end{equation}
As a consequence one has that for all $u \in H^2(M)$
\begin{equation}\label{eq:adaneg}
  \log \int_M e^{4(u - \ov{u})} dV_g \leq C + \frac{1}{8 \pi^2}
    \langle P_g^+ u, u \rangle.
\end{equation}
\end{lem}

\

\noindent From this result we derive an improved inequality for
functions which are concentrated at more than a single point,
related to a result in \cite{cl}. A consequence of this inequality
is that it allows us to give an upper bound (depending on $\int_M
Q_g dV_g$) for the number of concentration points of $e^{4u}$, where
$u$ is any function in $H^2(M)$ on which $II$ attains large negative
values, see Lemma \ref{l:II<-M}.

\begin{lem}\label{l:imprc}
For a fixed integer $\ell$, let $\Omega_1, \dots, \Omega_{\ell+1}$
be subsets of $M$ satisfying $dist(\Omega_i,\Omega_j) \geq \d _0$
for $i \neq j$, where $\d_0$ is a positive real number, and let
$\g_0 \in \left( 0, \frac{1}{\ell+1} \right)$. Then, for any
$\tilde{\e}
> 0$ and any $S > 0$ there exists a constant $C = C(\ell, \tilde{\e}, S,
\d_0, \g_0)$ such that
$$
  \log \int_M e^{4(u - \ov{u})} dV_g \leq C + \frac{1}{8 (\ell+1)
  \pi^2 - \tilde{\e}} \langle P_g u, u \rangle
$$
for all the functions $u \in H^2(M)$ satisfying
\begin{equation}\label{eq:ddmm}
    \frac{\int_{\Omega_i} e^{4u} dV_g}{\int_M e^{4u} dV_g} \geq \g_0,
  \quad \forall \; i \in \{1, \dots, \ell+1\}; \qquad \qquad \quad
  \sum _{i=1}^{\ov{k}} \a_i^2 \leq S.
\end{equation}
Here $\hat{u} = \sum_{i=1}^{\ov{k}} \a_i \hat{v}_i$ denotes the
component of $u$ in $V$.
\end{lem}

\begin{pf}
We modify the argument in \cite{cl} avoiding the use of truncations,
which is not allowed in the $H^2$ setting. Assuming without loss of
generality that $\ov{u} = 0$, we can find $\ell + 1$ functions $g_1,
\dots, g_{\ell+1}$ satisfying the following properties
\begin{equation}\label{eq:gi}
    \left\{%
\begin{array}{ll}
    g_i(x) \in [0,1] & \hbox{ for every } x \in M; \\
    g_i(x) = 1, & \hbox{ for every } x \in \Omega_i, i = 1, \dots, \ell+1; \\
    g_i(x) = 0, & \hbox{ if } dist(x,\Omega_i) \geq \frac{\d_0}{4}; \\
    \|g_i\|_{C^4(M)} \leq C_{\d_0}, &  \\
\end{array}%
\right.
\end{equation}
where $C_{\d_0}$ is a positive constant depending only on $\d_0$. By
interpolation, see \cite{Lio}, since $P_g^+$ is non-negative with
$ker P_g^+ = \R$, for any $\e
> 0$ there exists $C_{\e,\d_0}$ (depending only on $\e$ and $\d_0$)
such that, for any $v \in H^2(M)$ and for any $i \in \{1, \dots,
\ell + 1\}$ there holds
\begin{equation}\label{eq:pgiv}
    \langle P_g^+ g_i v, g_i v \rangle \leq \int_M g_i^2 (P_g^+ v, v)
    dV_g + \e \langle P_g^+ v, v \rangle + C_{\e,\d_0} \int_M v^2 dV_g.
\end{equation}
If we write $u$ as $u = u_1 + u_2$ with $u_1 \in L^\infty(M)$, then
from our assumptions we deduce
\begin{equation}\label{eq:inteu2}
    \int_{\Omega_i} e^{4 u_2} dV_g \geq e^{- 4 \|u_1\|_{L^\infty(M)}}
    \int_{\Omega_i} e^{4 u} dV_g \geq e^{- 4 \|u_1\|_{L^\infty(M)}}
    \g_0 \int_M e^{4 u} dV_g; \qquad i = 1, \dots, \ell + 1.
\end{equation}
Using \eqref{eq:gi}, \eqref{eq:inteu2} and then \eqref{eq:adaneg} we
obtain
\begin{eqnarray*}
    \log \int_M e^{4u} dV_g & \leq & \log \frac{1}{\g_0} + 4
    \|u_1\|_{L^\infty(M)} + \log \int_M e^{4 g_i u_2}
    dV_g + C \\ & \leq & \log \frac{1}{\g_0} + 4
    \|u_1\|_{L^\infty(M)} + C + \frac{1}{8 \pi^2}
    \langle P_g^+ g_i u_2, g_i u_2 \rangle + 4 \ov{g_i u_2},
\end{eqnarray*}
where $C$ depends only on $M$. We now choose $i$ such that $\langle
P_g^+ g_i u_2, g_i u_2 \rangle \leq \langle P_g^+ g_j u_2, g_j u_2
\rangle$ for every $j \in \{1, \dots, \ell + 1\}$. Since the
functions $g_1, \dots g_{\ell+1}$ have disjoint supports, the last
formula and \eqref{eq:pgiv} imply
$$
  \log \int_M e^{4u} dV_g \leq \log \frac{1}{\g_0} +  4
    \|u_1\|_{L^\infty(M)} + C + \left( \frac{1}{8(\ell+1)
  \pi^2} + \e \right) \langle P_g^+ u_2, u_2 \rangle + C_{\e,\d_0}
  \int_M u_2^2 dV_g + 4 \ov{g_i u_2}.
$$
Next we choose $\l_{\e,\d_0}$ to be an eigenvalue of $P_g^+$ such
that $\frac{C_{\e,\d_0}}{\l_{\e,\d_0}} < \e$, where $C_{\e,\d_0}$ is
given in the last formula, and we set
$$
  u_1 = P_{V_{\e,\d_0}} u; \qquad \qquad u_2 = P_{V_{\e,\d_0}^\perp} u,
$$
where $V_{\e,\d_0}$ is the direct sum of the eigenspaces of $P_g^+$
with eigenvalues less or equal to $\l_{\e,\d_0}$, and
$P_{V_{\e,\d_0}}, P_{V_{\e,\d_0}^\perp}$ denote the projections onto
$V_{\e,\d_0}$ and $V_{\e,\d_0}^\perp$ respectively. Since $\ov{u} =
0$, the $L^2$-norm and the $L^\infty$-norm on $V_{\e,\d_0}$ are
equivalent (with a proportionality factor which depends on $\e$ and
$\d_0$), and hence by our choice of $u_1$ and $u_2$ there holds
$$
  \|u_1\|_{L^\infty(M)}^2 \leq \hat{C}_{\e,\d_0} \langle P_g^+ u_1, u_1
  \rangle; \qquad \qquad C_{\e,\d_0} \int_M u_2^2 dV_g \leq
  \frac{C_{\e,\d_0}}{\l_{\e,\d_0}} \langle P_g^+ u_2, u_2 \rangle <
  \e \langle P_g^+ u_2, u_2 \rangle,
$$
where $\hat{C}_{\e,\d_0}$ depends on $\e$ and $\d_0$. Furthermore,
by the positivity of $P_g^+$ and the Poincar\'e inequality (recall
that $\ov{u} = 0$), we have
$$
  \ov{g_i u_2} \leq C \|u_2\|_{L^2(M)} \leq C \|u\|_{L^2(M)} \leq
  C \langle P_g^+ u, u \rangle^{\frac{1}{2}}.
$$
Hence the last formulas imply
\begin{eqnarray*}
  \log \int_M e^{4u} dV_g & \leq & \log \frac{1}{\g_0} +
  4 \hat{C}_{\e,\d_0} \langle P_g^+ u_1, u_1 \rangle^{\frac 12} + C + \left( \frac{1}{8(\ell+1)
  \pi^2} + \e \right) \langle P_g^+ u_2, u_2 \rangle + \e \langle P_g^+ u_2, u_2
  \rangle \\ & + & C \langle P_g^+ u_2, u_2 \rangle^{\frac 12} \leq \left(
  \frac{1}{8(\ell+1) \pi^2} + 3 \e \right) \langle P_g^+ u, u \rangle +
  \ov{C}_{\e,\d_0} + C + \log \frac{1}{\g_0},
\end{eqnarray*}
where $\ov{C}_{\e,\d_0}$ depends only on $\e$ and $\d_0$ (and
$\ell$, which is fixed). Now, since by \eqref{eq:ddmm} we have
uniform bounds on $\hat{u}$, we can replace $\langle P_g^+ u, u
\rangle$ by $\langle P_g u, u \rangle$ plus a constant in the
right-hand side. This concludes the proof.
\end{pf}

\

\noindent In the next lemma we show a criterion which implies the
situation described by the first condition in \eqref{eq:ddmm}.

\begin{lem}\label{l:er}
Let $\ell$ be a given positive integer, and suppose that $\e$ and
$r$ are positive numbers. Suppose that for a non-negative function
$f \in L^1(M)$ with $\|f\|_{L^1(M)} = 1$ there holds
$$
  \int_{\cup_{i=1}^\ell B_r(p_i)} f dV_g < 1 - \e \qquad \qquad \hbox{ for every
  $\ell$-tuple }
  p_1, \dots, p_\ell \in M.
$$
Then there exist $\ov{\e} > 0$ and $\ov{r} > 0$, depending only on
$\e, r, \ell$ and $M$ (but not on $f$), and $\ell + 1$ points
$\ov{p}_1, \dots, \ov{p}_{\ell+1} \in M$ (which depend on $f$)
satisfying
$$
  \int_{B_{\ov{r}}(\ov{p}_1)} f dV_g \geq \ov{\e}, \; \dots, \;
  \int_{B_{\ov{r}}(\ov{p}_{\ell+1})} f dV_g \geq \ov{\e}; \qquad \qquad
  B_{2 \ov{r}}(\ov{p}_i) \cap B_{2 \ov{r}}(\ov{p}_j) = \emptyset
  \hbox{ for } i \neq j.
$$
\end{lem}

\begin{pf}
Suppose by contradiction that for every $\ov{\e}, \ov{r} > 0$ there
is $f$ satisfying the assumptions and such that for every
$(\ell+1)$-tuple of points $p_1, \dots, p_{\ell+1}$ in $M$ we have
the implication
\begin{equation}\label{eq:sepk}
  \int_{B_{\ov{r}}(p_1)} f dV_g \geq \ov{\e}, \; \dots, \;
  \int_{B_{\ov{r}}(p_{\ell+1})} f dV_g \geq \ov{\e} \qquad
  \Rightarrow \qquad B_{2 \ov{r}}(p_i) \cap B_{2
  \ov{r}}(p_j) \neq \emptyset \hbox{ for some } i \neq j.
\end{equation}
We let $\ov{r} = \frac r8$, where $r$ is given in the statement. We
can find $h \in \N$ and $h$ points $x_1, \dots, x_h \in M$ such that
$M$ is covered by $\cup_{i=1}^h B_{\ov{r}}(x_i)$. For $\e$ given in
the statement of the Lemma, we also set $\ov{\e} = \frac{\e}{2 h}$.
We point out that the choice of $\ov{r}$ and $\ov{\e}$ depends on
$r, \e, \ell$ and $M$ only, as required.

Let $\{\tilde{x}_1, \dots, \tilde{x}_j \} \subseteq \{x_1, \dots,
x_h\}$ be the points for which $\int_{B_{\ov{r}}(\tilde{x}_i)} f
dV_g \geq \ov{\e}$. We define $\tilde{x}_{j_1} = \tilde{x}_1$, and
let $A_1$ denote the set
$$
  A_1 = \left\{ \cup_i B_{\ov{r}}(\tilde{x}_i) \; : \; B_{2 \ov{r}}(\tilde{x}_i)
  \cap B_{2 \ov{r}}(\tilde{x}_{j_1}) \neq \emptyset \right\}
  \subseteq B_{4 \ov{r}}(\tilde{x}_{j_1}).
$$
If there exists $\tilde{x}_{j_2}$ such that $B_{2
\ov{r}}(\tilde{x}_{j_2}) \cap B_{2 \ov{r}}(\tilde{x}_{j_1}) =
\emptyset$, we define
$$
  A_2 = \left\{ \cup_i B_{\ov{r}}(\tilde{x}_i) \; : \; B_{2 \ov{r}}(\tilde{x}_i)
  \cap B_{2 \ov{r}}(\tilde{x}_{j_2}) \neq \emptyset \right\}
  \subseteq B_{4 \ov{r}}(\tilde{x}_{j_2}).
$$
Proceeding in this way, we define recursively some points
$\tilde{x}_{j_3}, \tilde{x}_{j_4}, \dots, \tilde{x}_{j_s}$
satisfying
$$
  B_{2 \ov{r}}(\tilde{x}_{j_s}) \cap B_{2 \ov{r}}(\tilde{x}_{j_a})
  = \emptyset \; \qquad \forall \; 1 \leq a < s;
$$
and some sets $A_3, \dots, A_s$ by
$$
A_s = \left\{
  \cup_i B_{\ov{r}}(\tilde{x}_i) \; : \; B_{2 \ov{r}}(\tilde{x}_i)
  \cap B_{2 \ov{r}}(\tilde{x}_{j_s}) \neq \emptyset \right\}
  \subseteq B_{4 \ov{r}}(\tilde{x}_{j_s}).
$$
By \eqref{eq:sepk}, the process cannot go further than
$\tilde{x}_{j_\ell}$, and hence $s \leq \ell$. Using the definition
of $\ov{r}$ we obtain
\begin{equation}\label{eq:txitx1}
    \cup_{i=1}^j B_{\ov{r}}(\tilde{x}_i) \subseteq \cup_{i=1}^s A_i
    \subseteq \cup_{i=1}^s
    B_{4 \ov{r}} (\tilde{x}_{j_i}) \subseteq \cup_{i=1}^s
    B_{r} (\tilde{x}_{j_i}).
\end{equation}
Then by our choice of $h$, $\ov{\e}$, $\{\tilde{x}_1, \dots,
\tilde{x}_{j}\}$ and by \eqref{eq:txitx1} there holds
$$
  \int_{M \setminus \cup_{i=1}^s B_r(\tilde{x}_{j_i})} f dV_g \leq
  \int_{M \setminus \cup_{i=1}^j B_{\ov{r}}(\tilde{x}_i)} f dV_g
  \leq \int_{\left( \cup_{i=1}^h B_{\ov{r}} (x_i) \right) \setminus
  \left( \cup_{i=1}^j B_{\ov{r}}(\tilde{x}_i)
   \right)} f dV_g < (h-j) \ov{\e} \leq \frac{\e}{2}.
$$
Finally, if we chose $p_i = \tilde{x}_{j_i}$ for $i = 1, \dots, s$
and $p_i = \tilde{x}_{j_s}$ for $i = s+1, \dots, \ell$, we get a
contradiction to the assumptions of the lemma.
\end{pf}

\

\noindent Next we characterize some functions in $H^2(M)$ for which
the value of $II$ is large negative. Recall that the number $k$ is
given in formula \eqref{eq:kp3} and that $\hat{u}$ is the projection
of u on the direct sum of the eigenspaces of $P_g$ corresponding to
negative eigenvalues.

\begin{lem}\label{l:II<-M}
Under the assumptions of Theorem \ref{th:ex}, and for $k_P \in (8 k
\pi^2, 8(k+1) \pi^2)$ with $k \geq 1$, the following property holds.
For any $S > 0$, any  $\e > 0$ and any $r > 0$ there exists a large
positive $L = L(S, \e, r)$ such that for every $u \in H^2(M)$ with
$II(u) \leq - L$ and $\|\hat{u}\| \leq S$ there exists $k$ points
$p_{1,u}, \dots, p_{k,u} \in M$ such that
\begin{equation}\label{eq:caz}
 \int_{M \setminus \cup_{i=1}^{k} B_r(p_{i,u})} e^{4u} dV_g < \e.
\end{equation}
\end{lem}

\begin{pf}
Suppose by contradiction that the statement is not true, namely that
there exist $S, \e, r > 0$ and $(u_n)_n \subseteq H^2(M)$ with
$\|\hat{u}_n\| \leq S$, $II(u_n) \to - \infty$ and such that for
every $k$-tuple $p_1, \dots, p_k$ in $M$ there holds
$\int_{\cup_{i=1}^k B_r(p_i)} e^{4 u_n} dV_g < 1 - \e$. Recall that
without loss of generality, since $II$ is invariant under
translation by constants in the argument, we can assume that for
every $n$ $\int_M e^{4 u_n} dV_g = 1$. Then we can apply Lemma
\ref{l:er} with $\ell = k$, $f = e^{4u_n}$, and in turn Lemma
\ref{l:imprc} with $\d_0 = 2 \ov{r}$, $\Omega_1 =
B_{\ov{r}}(\ov{p}_1)$, \dots, $\Omega_{k+1} =
B_{\ov{r}}(\ov{p}_{k+1})$ and $\g_0 = \ov{\e}$, where $\ov{\e}$,
$\ov{r}$ and $(\ov{p}_i)_i$ are given by Lemma \ref{l:er}. This
implies that for any given $\tilde{\e}
> 0$ there exists $C > 0$ depending only on $S, \e, \tilde{\e}$ and
$r$ such that
\begin{eqnarray*}
    II(u_n) \geq \langle P_g u_n, u_n \rangle + 4 \int_M Q_g u_n dV_g - C k_P -
    \frac{k_P}{8 (k+1) \pi^2 - \tilde{\e}} \langle P_g u_n, u_n \rangle - 4 k_P
    \ov{u}_n,
\end{eqnarray*}
where $C$ is independent of $n$. Since $k_P < 8 (k+1) \pi^2$, we can
choose $\tilde{\e} > 0$ so small that $1 - \frac{k_P}{8 (k+1) \pi^2
- \tilde{\e}} := \d > 0$. Hence using also the Poincar\'e inequality
we deduce
\begin{eqnarray}{\label{eq:tech}}
  II(u_n) & \geq & \d \langle P_g u_n, u_n \rangle + 4 \int_M Q_g (u_n - \ov{u}_n)
  dV_g - C k_P \nonumber \\ & \geq & \d \langle P_g u_n, u_n \rangle - 4 C
  \langle P_g u_n, u_n \rangle^{\frac{1}{2}} - C k_P \geq - C.
\end{eqnarray}
This violates our contradiction assumption, and concludes the proof.
\end{pf}

\section{Mapping sublevels of $II$ into $A_{k,\ov{k}}$}\label{s:msub}

In this section we show how to map non trivially some sublevels of
the functional $II$ into the set $A_{k,\ov{k}}$. Since adding a
constant to the argument of $II$ does not affect its value, we can
always assume that the functions $u \in H^2(M)$ we are dealing
with satisfy the normalization \eqref{eq:noul} (with $u$ instead
of $u_l$). Our goal is to prove the following result.

\begin{pro}\label{p:map}
For $k \geq 1$ (see \eqref{eq:kp3}) there exists a large $L > 0$
and a continuous map $\Psi$ from the sublevel $\{ II < - L \}$
into $A_{k,\ov{k}}$ which is topologically non-trivial. For $k_P <
8 \pi^2$ and $\ov{k} \geq 1$ the same is true with $A_{k,\ov{k}}$
replaced by $S^{\ov{k}-1}$
\end{pro}

\noindent We divide the section into two parts. First we derive some
properties of the set $M_k$ for $k \geq 1$. Then we turn to the
construction of the map $\Psi$. Its non-triviality  will follow from
Proposition \ref{p:mkm} below, where we show that there is another
map $\Phi$ from $A_{k,\ov{k}}$ into $H^2(M)$ such that $\Psi \circ
\Phi$ is homotopic to the identity on $A_{k,\ov{k}}$, which is not
contractible by Corollary \ref{c:akk}.

\subsection{Some properties of the set $M_k$}\label{ss:prop}

In this subsection we collect some useful properties of the set
$M_k$, beginning with some local ones near the singularities, namely
the subsets $M_j \subseteq M_k$ with $j < k$. Although the
topological structure of the barycenters is well-known, we need some
estimates of quantitative type concerning the metric distance. The
reason, as mentioned in the introduction, is that the amount of
concentration of $e^{4u}$ (where $u \in \{ II \leq - L \}$, see
Lemma \ref{l:II<-M}) near a single point can be arbitrarily small.
In this way we are forced to define a projection which depends on
{\em all} the distances from the $M_j$'s, see Subsection
\ref{ss:cinstr}, which requires some preliminary considerations. We
recall that on $M_k$ we are adopting the metric induced by
$C^1(M)^*$, see Section \ref{s:pr}, and for $j < k$ we set $d_j(\s)
= dist(\s,M_j)$, $\s \in M_k$. Then for $\e > 0$ and $2 \leq j \leq
k$, we define
$$
  M_j(\e) = \left\{ \s \in M_j \; : \; d_{j-1}(\s) > \e
  \right\}.
$$
For convenience, we extend the definition also to the case $j =
1$, setting
$$
  M_1(\e) : = M_1.
$$
We give a first quantitative description of the set $M_j(\e)$,
which leads immediately to (the known) Corollary \ref{cor:Mksmo}.

\begin{lem}\label{l:Pj-1}
Let $j \in \{2, \dots, k\}$. Then there exists $\e$ sufficiently
small with the following property. If $\s \in M_j(\e)$, $\s =
\sum_{i=1}^j t_i \d_{x_i}$, then there holds
\begin{equation}\label{eq:j-1j-1}
  t_i \geq \frac{\e}{2}; \qquad dist(x_i, x_l) \geq \frac{\e}{2};
  \qquad \qquad i, l = 1, \dots, j, i \neq l.
\end{equation}
\end{lem}

\begin{pf}
Let $\s = \sum_{i=1}^j t_i \d_{x_i} \in M_j(\e)$. Assuming by
contradiction that the first inequality in \eqref{eq:j-1j-1} is not
satisfied, there would exists $\ov{\i} \in \{1, \dots, j\}$ such
that $t_{\ov{\i}} < \frac \e 2$. Then, for $\tilde{\i} \in \{1,
\dots, j\}$, $\tilde{\i} \neq \ov{\i}$, we consider the following
element
$$
  \hat{\s} = (t_{\ov{\i}} + t_{\tilde{\i}}) \d_{x_{\tilde{\i}}} +
  \sum_{i = 1, \dots, j, \; i \neq \ov{\i}, \tilde{\i}} t_i \d_{x_i} \in M_{j-1}.
$$
For any function $f \in C^1(M)$ with $\|f\|_{C^1(M)} \leq 1$ there
holds clearly
$$
  \left| (\s, f) - (\hat{\s}, f) \right| \leq t_{\ov{\i}} \left(
  |f(x_{\ov{\i}})| + |f(x_{\tilde{\i}})| \right) \leq 2 t_{\ov{\i}}.
$$
Taking the supremum with respect to $f$ we deduce
$$
  \e < dist(\s, M_{j-1}) \leq dist(\s, \hat{\s}) = \sup_f \left|
  (\s, f) - (\hat{\s}, f) \right| \leq 2 t_{\ov{\i}}.
$$
This gives us a contradiction. Let us prove now the second
inequality. Assuming that there are $x_i, x_l \in M$ with, $x_i \neq
x_l$ and $dist(x_i, x_l) < \frac \e 2$ (for $\e$ sufficiently
small), let us define the element
$$
  \hat{\s} = (t_{i} + t_{l}) \d_{\frac 12 x_{i}
  + \frac 12 x_l} + \sum_{s = 1, \dots, j, \; s \neq i, l}
  t_s \d_{x_s} \in M_{j-1},
$$
see the notation introduced in Section \ref{s:pr} for the convex
combination of the points $x_i$ and $x_l$. Similarly as before, for
$\|f\|_{C^1(M)} \leq 1$ we obtain
\begin{eqnarray*}
    \left| (\s, f) - (\hat{\s}, f) \right| & \leq & t_{i}
  \left|f(x_{i}) - f\left(\frac{x_{i} + x_{l}}{2}
  \right) \right| + t_{l} \left|f(x_{l}) - f
  \left(\frac{x_{i}+x_{l}}{2} \right) \right| \\
  & \leq & \left|f(x_{i}) - f\left(\frac{x_{i} +
  x_{l}}{2} \right) \right| + \left|f(x_{l})
  - f \left(\frac{x_{i} + x_{l}}{2} \right) \right|.
\end{eqnarray*}
Taking the supremum with respect to such functions $f$, since they
all have Lipschitz constant less or equal than 1, we deduce
$$
  \e < dist(\s, M_{j-1}) \leq dist(\s, \hat{\s}) = \sup_f \left|
  (\s, f) - (\hat{\s}, f) \right| \leq 2 dist(x_{i}, x_{l}).
$$
This gives us a contradiction and concludes the proof.
\end{pf}

\begin{cor}\label{cor:Mksmo} (well-known)
The set $M_1$ is a smooth manifold in $C^1(M)^*$. Furthermore, for
any $\e > 0$ and for $j \geq 2$, the set $M_j(\e)$ is also a smooth
(open) manifold of dimension $5j - 1$.
\end{cor}

\begin{pf}
The first assertion is obvious. Regarding the second one, the
previous lemma guarantees that all the numbers $t_i$ are uniformly
bounded away from zero and that the mutual distance between the
points $x_i$'s is also uniformly bounded from below. Therefore,
recalling that the $t_i$'s satisfy the constraint $\sum_i t_i =
1$, each element of $M_j(\e)$ can be smoothly parameterized by $4
j$ coordinates locating the points $x_i$'s and by $j - 1$
coordinates identifying the numbers $t_i$'s.
\end{pf}

\

\noindent We show next that it is possible to define a continuous
homotopy which brings points in $M_k$, which are close to $M_j(\e)$,
onto $M_j\left( \frac \e2 \right)$. We also provide some
quantitative estimates on the deformation. Our goal is to patch
together projections onto sets of different dimensions (of the form
$M_j(\e_j)$ for suitable $\e_j$'s), as shown in Figure \ref{fig2}
below. The proof of Lemma \ref{l:Pjl} is postponed to the appendix.

\begin{lem}\label{l:Pjl}
Let $j \in \{1, \dots, k-1\}$, and let $\e > 0$. Then there exist
$\hat{\e}$ ($\ll \e^2$), depending only on $\e$ and $k$, and a map
$T^t_j$, $t \in [0,1]$, from the set
$$
\hat{M}_{k,j}^{\hat{\e},\e} := \left\{ \s \in M_k \; : \; dist(\s,
M_j(\e)) < \hat{\e} \right\}
$$
into $M_k$ such that the following five properties hold true
\begin{description}
    \item[(i)] $T^0_j = Id$ and $T^t_j|_{M_j} = Id|_{M_j}$ for every $t \in [0,1]$;
    \item[(ii)] $T^1_j(\s) \in M_j \left( \frac{\e}{2}
    \right)$ for every $\s \in \hat{M}_{k,j}^{\hat{\e},\e}$;
    \item[(iii)] $dist(T^0_j(\s), T^t_j(\s)) \leq
    C_{k,\e} \sqrt{\hat{\e}}$ for every $\s \in \hat{M}_{k,j}^{\hat{\e},\e}$ and
    for every $t \in [0,1]$;
    \item[(iv)] if $\s \in \hat{M}_{k,j}^{\hat{\e},\e} \cap M_l$ for
    some $l \in \{j, \dots, k-1\}$, then $T^t_j(\s) \in M_l$ for every
    $t \in [0,1]$;
    \item[(v)] if $\s \in \hat{M}_{k,j}^{\hat{\e},\e} \cap M_j$, then
    $T^t_j(\s) = \s$ for every $t \in [0,1]$.
\end{description}
The constant $C_{k,\e}$ in {\bf (iii)} depends only on $k$ and $\e$,
and not on $t$ and $\hat{\e}$ (provided the latter is sufficiently
small).
\end{lem}

\begin{rem}\label{r:samej}
We notice that, by the property {\bf (iv)} in the statement of Lemma
\ref{l:Pjl}, the above homotopy is well defined also from each $M_l$
into itself, for $l \in \{1, \dots, k-1\}$, and extends continuously
to a neighborhood of $M_l$ in $M_k$.
\end{rem}

\

\noindent Since $M_j\left( \frac{\e}{4} \right)$ is a smooth
finite-dimensional manifold in $C^1(M)^*$ by Corollary
\ref{cor:Mksmo}, we can define a continuous projection $P_j$ (see
the comments at the beginning of the proof of Lemma \ref{l:Pjl} in
the appendix) from $\hat{M}_{k,j}^{\hat{\e},\e}$ into $M_j\left(
\frac{\e}{2} \right)$, which is compactly contained in $M_j\left(
\frac{\e}{4} \right)$. We have then an immediate consequence of the
previous lemma.

\begin{cor}\label{c:hom}
Let $T^t_{j}$ denote the map constructed in Lemma \ref{l:Pjl} above.
Then for $\hat{\e}$ sufficiently small there exists an homotopy
$H_j^t$, $t \in [0,1]$, between $T^1_{j}(\s)$ and $P_j(\s)$ within
$M_j \left( \frac{\e}{2} \right)$, namely a map satisfying the
following properties
\begin{equation}\label{eq:Hjt}
    \left\{%
\begin{array}{ll}
    H_j^t(\s) \in M_j \left( \frac{\e}{2} \right) & \hbox{ for every }
    t \in [0,1] \hbox{ and every } \s \in \hat{M}^{\hat{\e},\e}_{k,j}; \\
    H_j^0(\s) = T^1_j(\s) & \hbox{ for every } \s \in \hat{M}^{\hat{\e},\e}_{k,j}; \\
    H_j^1(\s) = P_j(\s) & \hbox{ for every } \s \in \hat{M}^{\hat{\e},\e}_{k,j}. \\
\end{array}%
\right.
\end{equation}
\end{cor}


\

\noindent In view of Corollary \ref{c:hom}, we can also modify the
map $T^t_j$ by composing it with the above homotopy $H^t_j$, namely
we set
\begin{equation}\label{eq:hatT}
    \hat{T}^t_j(\s) = \left\{%
\begin{array}{ll}
    T^{2t}_j(\s), & \hbox{ for } t \in \left[ 0, \frac 12 \right]; \\
    H^{2t-1}_j \circ T^1_j(\s), & \hbox{ for } t \in \left[ \frac 12, 1 \right].  \\
\end{array}%
\right.
\end{equation}
In this way, if for some element of $C^1(M)^*$ both the projections
$P_j$ and $P_l$ are defined, for $1 \leq j < l \leq k$, composing
$\hat{T}^t_j$ with $P_l$ we obtain an homotopy between $P_l$ and
$P_j$ within $M_l \cap M_j\left( \frac \e 2 \right)$, see Remark
\ref{r:samej}. This fact will be crucially used in the proof of
Lemma \ref{l:defmap} below.

\

\noindent Next we recall the following result, which is necessary
in order to carry out the topological argument below. For
completeness, we give a brief idea of the proof.

\begin{lem}\label{l:nonc} (well-known)
For any $k \geq 1$, the set $M_k$ is non-contractible.
\end{lem}

\begin{pf}
For $k = 1$ the statement is obvious, so we consider the case $k
\geq 2$. The set $M_k \setminus M_{k-1}$, see Corollary
\ref{cor:Mksmo}, is an open manifold of dimension $5k - 1$. It is
possible to prove that, even if $M_{k-1}$ is not a smooth manifold
(for $k \geq 3)$, it is anyway a Euclidean Neighborhood Retract,
namely it is a contraction of some of its neighborhoods which has
smooth boundary (of dimension $5k - 2$), see \cite{bc},
\cite{bred}. Therefore $M_k$ has an orientation (mod 2) with
respect to $M_{k-1}$, namely the relative homology class
$H_{5k-1}(M_k, M_{k-1}; \mathbb{Z}_2)$ is non-trivial. Consider
now this part of the exact homology sequence of the pair $(M_k,
M_{k-1})$
$$
  \cdots \to H_{5k-1}(M_{k-1};\mathbb{Z}_2) \to H_{5k-1}(M_k;\mathbb{Z}_2)
  \to H_{5k-1}(M_k, M_{k-1};\mathbb{Z}_2) \to
  H_{5k-2}(M_{k-1};\mathbb{Z}_2) \to \cdots
$$
Since the dimension of (the stratified set) $M_{k-1}$ is less or
equal than $5(k-1)-1 < 5k - 2$, both the homology groups
$H_{5k-1}(M_{k-1};\mathbb{Z}_2)$ and
$H_{5k-2}(M_{k-1};\mathbb{Z}_2)$ vanish, and therefore
$H_{5k-1}(M_k;\mathbb{Z}_2) \simeq H_{5k-1}(M_k,
M_{k-1};\mathbb{Z}_2) \neq 0$. The proof is concluded.
\end{pf}

\

\noindent From the preceding Lemma and from a standard application
of the Majer-Vietoris Theorem is it easy to deduce the following
result.

\begin{cor}\label{c:akk}
For any (relative) integers $k \geq 1$ and $\ov{k} \geq 0$, the
set $A_{k,\ov{k}}$ is non-contractible.
\end{cor}

\subsection{Construction of $\Psi$}\label{ss:cinstr}

In this subsection we finally construct the map $\Psi$, using the
preceding results about the set $M_k$. First we show how to
construct some partial projections on the sets $M_j(\e)$ for $\e >
0$. When referring to the distance of a function in $L^1(M)$ from a
set $M_j$, we always adopt the metric induced by $C^1(M)^*$. The
comments before Corollary \ref{c:hom} yield the following result.

\begin{lem}\label{l:Pj}
Suppose that $f \in L^1(M)$, $f \geq 0$ and that $\int_M f dV_g =
1$. Then, given any $\e > 0$ and any $j \in \{1, \dots, k\}$,
there exists $\hat{\e}
> 0$, depending on $j$ and $\e$ with the following property. If
$dist(f, M_{j}(\e)) \leq \hat{\e}$, then there is a continuous
projection $P_j$ mapping $f$ onto $M_j \left( \frac{\e}{2} \right)$.
\end{lem}


\

\noindent Next we define an auxiliary map $\hat{\Psi}$ from a
suitable sublevel of $II$ into $M_k$.

\begin{lem}\label{l:defmap}
For $k \geq 1$ there exist a large $\hat{L} > 0$ and a continuous
map $\hat{\Psi}$ from $\{ II \le - \hat{L} \} \cap \{ \|\hat{u}\|
\leq 1 \}$ into $M_k$. Here, as before, $\hat{u}$ denotes the
component of $u$ belonging to $V$, the direct sum of the negative
eigenspaces of $P_g$ (if any, otherwise we impose no restriction on
$u$ except for $II(u) \leq - \hat{L}$).
\end{lem}

\begin{pf}
First we define some numbers
$$
  \e_k \ll \e_{k-1} \ll \dots \ll \e_2 \ll \e_1 \ll 1
$$
in the following way. We choose $\e_1$ so small that there is a
projection $P_1$ from the non-negative $L^1(M)$ functions in an
$\e_1$-neighborhood of $M_1$ onto $M_1$ (by Lemma \ref{l:Pj}). We
now can apply again Lemma \ref{l:Pj} with $j = 2$, $\e = 4 \e_1$
and, obtaining the corresponding (sufficiently small) $\hat{\e}$, we
define $\e_2 = \frac{\hat{\e}}{4}$. Then we choose the numbers
$\e_3, \dots, \e_k$ iteratively in the same way.

For any $i = 1, \dots, k$, let $f_i$ be a smooth non-increasing
cutoff function which satisfies the following properties
\begin{equation}\label{eq:fi}
    \left\{%
\begin{array}{ll}
    f_i(t) = 1 & \hbox{ for } t \leq \e_i; \\
    f_i(t) = 0 & \hbox{ for } t \geq 2 \e_i. \\
\end{array}%
\right.
\end{equation}
Next we choose suitably the large number $\hat{L}$. In order to do
this, we apply Lemma \ref{l:II<-M} with $S = 1$ and some small $\e$.
It is easy to see that if $\e$ is chosen first sufficiently small,
and then $\hat{L} = L$ sufficiently large, then for any $u \in
H^2(M)$ with $II(u) \leq - \hat{L}$ (and $\int_M e^{4u} dV_g = 1$)
there holds $dist(e^{4u}, M_k) < \e_k$.

Now, given $u \in H^2(M)$ with $II(u) \leq - \hat{L}$, we let $j$
(depending on $u$) denote the first integer such that
$f_j(dist(e^{4u}, M_j)) = 1$. We notice that for $j > 1$, since
$f_{j-1}(dist(e^{4u}, M_{j-1})) < 1$, there holds $dist(e^{4u},
M_{j-1}) > \e_{j-1}$ and $dist(e^{4u}, M_{j-1}) < \e_j$.
Therefore, by Lemma \ref{l:Pj} and our choice of the $\e_i$'s, the
projection $P_j(e^{4u})$ is well-defined. Then we set
$$
  \hat{\Psi} (u) = \hat{T}_1^{f_1(dist(e^{4u},M_1))} \circ
  \hat{T}_2^{f_2(dist(e^{4u},M_2)} \circ \dots \circ
  \hat{T}_{j-1}^{f_{j-1}(dist(e^{4u},M_{j-1}))} \circ P_j(e^{4u}).
$$
Some comments are in order. This definition depends in principle on
the index $j$, which is a function of $u$. Nevertheless, since all
the distance functions from the $M_l$'s are continuous, and since
$\hat{T}^1_l = P_l$, see \eqref{eq:hatT}, the above map $\hat{\Psi}$
is indeed well defined and continuous in $u$, see Remark
\ref{r:samej} and the comments after Corollary \ref{c:hom}.
\end{pf}

\

\begin{figure}[ht]
\begin{center}
\includegraphics[width=15cm]{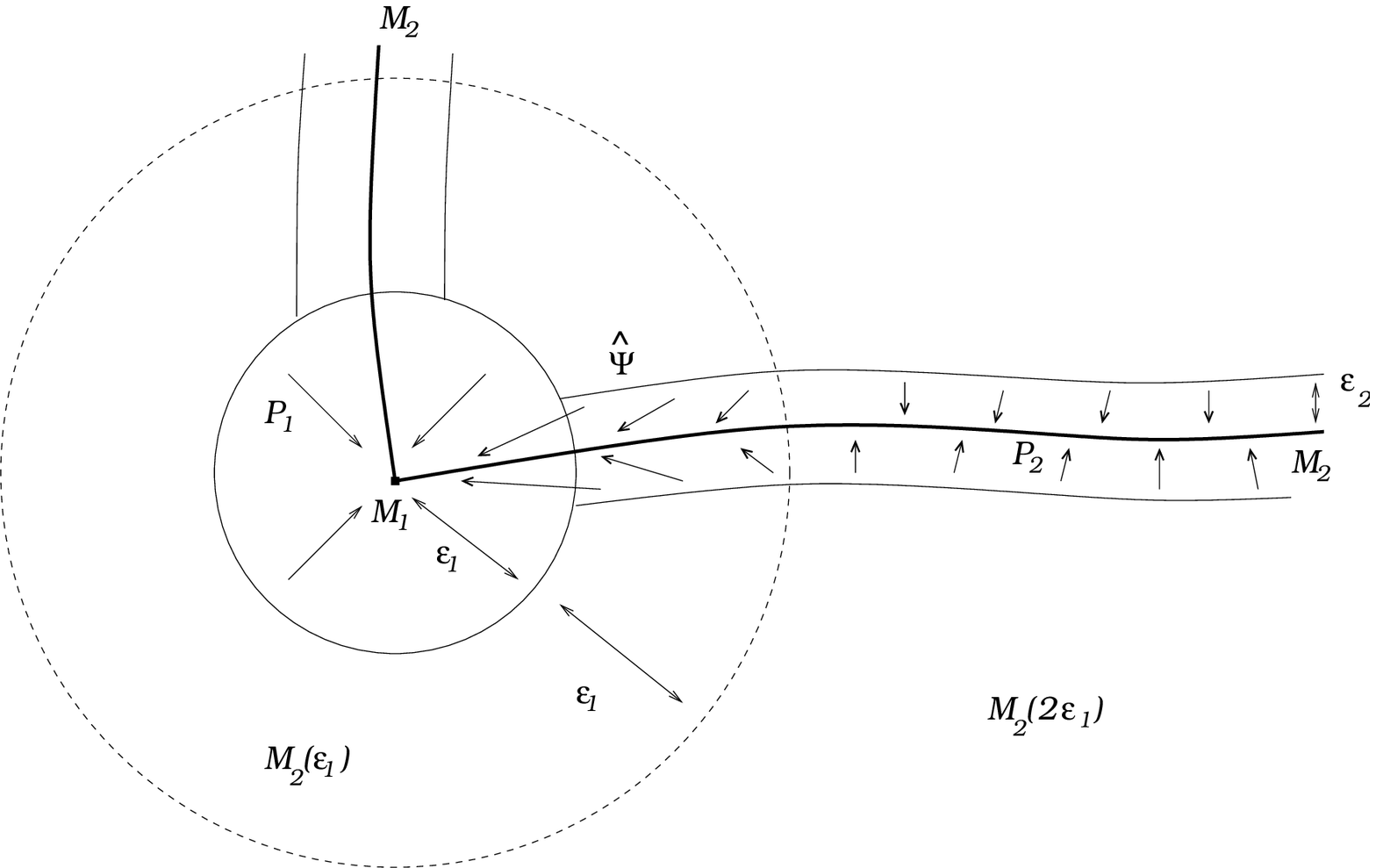}
\caption{the map $\hat{\Psi}$ for $k = 2$} \label{fig2}
\end{center}
\end{figure}

\noindent In Figure \ref{fig2} below we sketch the construction of
the map $\hat{\Psi}$ for the case $k = 2$, which is the simplest
among the non-trivial ones. $M_1$ is depicted as a single point,
while $M_2$ as a couple of curves. The region between the two
circles, which represents $\{ \e_1 \leq dist(e^{4u}, M_1) \leq 2
\e_1 \}$, is where the homotopy $\hat{T}^t_1$ (and hence the
construction of Lemma \ref{l:Pjl}) is used.

\

\noindent We are finally in position to introduce the global map
$\Psi$. If $\hat{v}_1, \dots, \hat{v}_{\ov{k}}$ form an orthonormal
basis (in $L^2(M)$) of $V$, $V$ being the direct sum of the
eigenspaces of $P_g$ corresponding to negative eigenvalues, see
Section \ref{s:pr}, we define the ${\ov{k}}$-vector
$$
  s(u) = \left( (\hat{v}_1, u)_{L^2(M)}, \dots, (\hat{v}_{\ov{k}},
  u)_{L^2(M)} \right) \in \R^{\ov{k}}.
$$
Then, if $\hat{L}$ is as in Lemma \ref{l:defmap} and if $\ov{\s}$ is
any fixed element of $M_k$, in the case $k \geq 1$ we let $\Psi : \{
II \leq - \hat{L} \} \to A_{k,\ov{k}}$ be defined by
\begin{equation}\label{eq:lk}
    \Psi(u) = \left\{%
\begin{array}{ll}
    (\hat{\Psi}(u), s(u)) & \hbox{ for } |s(u)| \leq 1; \\
    \left( \ov{\s}, \frac{s(u)}{|s(u)|} \right), & \hbox{ for }
    |s(u)| > 1. \\
\end{array}%
\right.
\end{equation}
Since for $|s|$ tending to $1$ the set $M_k$ is collapsing to a
single point in $A_{k,\ov{k}}$, see \eqref{eq:akovk}, the map $\Psi$
is continuous.

On the other hand, if $k_P < 8 \pi^2$ and if $\ov{k} \geq 1$ we
just set
\begin{equation}\label{eq:lk2}
  \Psi(u) = \frac{s(u)}{|s(u)|}.
\end{equation}

\

\begin{pfn} {\sc of Proposition \ref{p:map}.}
It remains only to prove the non-triviality of the map $\Psi$. This
follows from Corollary \ref{c:akk} and from {\bf (b)} in Proposition
\ref{p:mkm}.
\end{pfn}

\section{Mapping $A_{k,\ov{k}}$ into low sublevels of $II$}\label{s:test}

\noindent The next step consists in finding a map $\Phi$ from
$A_{k,\ov{k}}$ (resp. from $S^{\ov{k}-1}$) into $H^2(M)$ on which
image the functional $II$ attains large negative values.

\begin{pro}\label{p:mkm}
Let $\Psi$ be the map defined in the previous section. Then,
assuming $k \geq 1$ (resp. $k_P < 8 \pi^2$ and $\ov{k} \geq 1$), for
any $L > 0$ sufficiently large (such that Proposition \ref{p:map}
applies) there exists a map $\Phi_{\ov{S},\ov{\l}} : A_{k,\ov{k}}
\to H^2(M)$ (resp. $\Phi_{\ov{S}} : S^{\ov{k}-1} \to H^2(M)$) with
the following properties
\begin{description}
    \item[(a)] $II(\Phi_{\ov{S},\ov{\l}}(z)) \leq - L$ for any
    $z \in A_{k,\ov{k}}$ (resp. $II(\Phi_{\ov{S}}(z)) \leq - L$ for any
    $z \in S^{\ov{k}-1}$);
    \item[(b)] $\Psi \circ \Phi_{\ov{S},\ov{\l}}$ is homotopic to the
    identity on $A_{k,\ov{k}}$ (resp. $\Psi \circ \Phi_{\ov{S}}$ is
    homotopic to the identity on $S^{\ov{k}-1}$).
\end{description}
\end{pro}

\

\noindent In order to prove this proposition we need some
preliminary notations and lemmas. For $\d > 0$ small, consider a
smooth non-decreasing cut-off function $\chi_\d : \R_+ \to \R$
satisfying the following properties
\begin{equation}\label{eq:chid}
    \left\{%
\begin{array}{ll}
    \chi_\d(t) =  t & \hbox{ for } t \in [0,\d]; \\
    \chi_\d(t) =  2 \d & \hbox{ for } t \geq 2 \d; \\
    \chi_\d(t) \in [\d, 2 \d] & \hbox{ for } t \in [\d, 2 \d]. \\
\end{array}%
\right.
\end{equation}
Then, given $\s \in M_k$ $\left( \s = \sum_{i=1}^k t_i \d_{x_i}
\right)$ and $\l
> 0$, we define the function $\var_{\l,\s} : M \to \R$ as
\begin{equation}\label{eq:pls}
  \var_{\l,\s} (y) = \frac 14 \log
  \sum_{i=1}^k t_i \left( \frac{2 \l}{1 + \l^2 \chi_\d^2
  \left( d_i(y) \right)} \right)^4; \qquad \quad y \in M,
\end{equation}
where we have set
$$
 d_i(y) = dist(y,x_i), \qquad \qquad y \in M,
$$
with $dist(\cdot, \cdot)$ denoting the distance function on $M$.
We are now in position to define the function
$\Phi_{\ov{S},\ov{\l}} : A_{k,\ov{k}} \to H^2(M)$. For large
$\ov{S}$ and $\ov{\l}$ we let
$$
  \Phi_{\ov{S},\ov{\l}}(\s,s) =
  \left\{%
\begin{array}{ll}
    \var_{s} + \var_{\ov{\l},\s} & \hbox{ for } |s| \leq \frac 14;
    \\ \var_{s} + \var_{2 \ov{\l} - 1 + 4 (1 - \ov{\l}) |s|,\s} &
    \hbox{ for } \frac 14 \leq |s| \leq \frac 12; \\ \var_{s} + 2
    (1 - \var_{1,\s}) |s| + 2 \var_{1,\s} - 1 & \hbox{ for } |s|
    \geq \frac 12,  \\
\end{array}%
\right.
$$
where
$$
  s = (s_1, \dots, s_{\ov{k}}); \qquad \qquad \var_{s} (y) =
  \ov{S} \sum_{i=1}^{\ov{k}} s_i \hat{v}_i(y).
$$
For $k_P < 8 \pi^2$ and for $\ov{k} \geq 1$ we just set
$$
  \Phi_{\ov{S}}(s) = \var_s, \qquad \qquad |s| = 1.
$$

\noindent Notice that the map is well defined on $A_{k,\ov{k}}$.

\

\noindent We state now two preliminary lemmas, postponing the proof
of the first to the appendix.

\begin{lem}\label{l:pff}
Suppose $\var_{\l,\s}$ is as in \eqref{eq:pls}. Then as $\l \to +
\infty$ one has
$$
  \langle P_g \var_{\l,\s}, \var_{\l,\s} \rangle \leq
   \left( 32 k \pi^2 + o_\d(1) \right) \log \l + C_\d \qquad \quad
   (\hbox{uniformly in } \s \in M_k),
$$
where $o_\d(1) \to 0$ as $\d \to 0$, and where $C_\d$ is a constant
independent of $\l$ and $(x_i)_i$.
\end{lem}

\

\begin{lem}\label{l:foth}
For $k \geq 1$ (resp. for $k_P < 8 \pi^2$ and for $\ov{k} \geq 1$),
given any $L > 0$, there exist a small $\d$, some large $\ov{S}$ and
$\ov{\l}$ such that $II(\Phi_{\ov{S},\ov{\l}}(\s,s)) \leq - L$ for
every $(\s,s) \in A_{k,\ov{k}}$ (resp. $II(\Phi_{\ov{S}}(s)) \leq -
L$ for every $s \in S^{\ov{k}-1}$).
\end{lem}

\begin{pf}
We begin with the case $k \geq 1$, and we prove first the following
three estimates (recall that $\l_{\ov{k}}$ is the biggest negative
eigenvalue of $P_g$)
\begin{equation}\label{eq:estQvar}
    \int_M Q_g (\var_s + \var_{\l,\s}) dV_g = - k_P \log \l +
    O(\d^4 \log \l) + O(|\log \d|) + \ov{S} O(|s|) + O(1);
\end{equation}
\begin{equation}\label{eq:estexpvar}
    \log \int_M \exp \left( 4 (\var_s + \var_{\l,\s}) \right)
    dV_g = O(1) + O(\ov{S} |s|);
\end{equation}
\begin{equation}\label{eq:estPvar}
    \langle P_g (\var_s + \var_{\l,\s}),
    (\var_s + \var_{\l,\s}) \rangle \leq - |\l_{\ov{k}}| |s|^2 \ov{S}^2 +
    32 k \pi^2 (1 + o_\d(1)) \log \l + C_\d + O(\d^4 |s| \ov{S}).
\end{equation}

\

\noindent \underline{Proof of \eqref{eq:estQvar}}. We have
$$
  \var_{\l,\s} (y) = \log \frac{2 \l}{1 + 4 \l^2 \d^2},
  \quad \hbox{ for } y \in M \setminus \cup_{i=1}^k B_{2 \d} (x_i),
$$
and
$$
   \log \frac{2 \l}{1 + 4 \l^2 \d^2} \leq \var_{\l,\s}
   (y) \leq  \log 2 \l, \quad \hbox{ for } y \in
   \cup_{i=1}^k B_{2 \d} (x_i).
$$
Writing
\begin{eqnarray*}
    \int_M Q_g(y) \var_{\l,\s}(y) d V_g(y) & = & \log
    \frac{2 \l}{1 + 4 \l^2 \d^2} \int_M Q_g (y) dV_g (y) \\ & + &
    \int_M Q_g(y) \left( \var_{\l,\s}(y) - \log \frac{2
    \l}{1 + 4 \l^2 \d^2} \right) dV_g(y),
\end{eqnarray*}
from the last three formulas it follows that
\begin{equation}\label{eq:intpls}
    \int_M Q_g(y) \var_{\l,\s}(y) d V_g(y) =  k_P
    \log \frac{2 \l}{1 + 4 \l^2 \d^2} + O \left( \d^4 \log
    (1 + 4 \l^2 \d^2) \right).
\end{equation}
Furthermore recalling that the average of $\var_s$ is zero (since
all the $\hat{v}_i$'s have zero average, see Section \ref{s:pr}), we
also deduce that
\begin{equation}\label{eq:intps}
    \int_M Q_g(y) \var_s (y) dV_g(y) = \ov{S} \sum_{i=1}^{\ov{k}}
  s_i \int_M Q_g(y) \hat{v}_i (y) dV_g(y) = \ov{S} O(|s|).
\end{equation}
Hence \eqref{eq:intpls} and \eqref{eq:intps} yield
$$
  \int_M Q_g(y) (\var_s + \var_{\l,\s}(y)) dV_g(y) = k_P \log
  \frac{2 \l}{1 + 4 \l^2 \d^2} + O \left( \d^4 \log (1 + 4 \l^2 \d^2)
  \right) + \ov{S} O(|s|),
$$
which implies immediately \eqref{eq:estQvar}.

\

\noindent \underline{Proof of \eqref{eq:estexpvar}}. We recall that
in $V$ the $L^2$-norm and the $L^\infty$ norm are equivalent.
Therefore, noticing that
\begin{equation}\label{eq:xx}
    \exp \left( 4 (\var_s(y)) \right) \in \left[ \exp
  (4 \inf_M \var_s), \exp (4 \sup_M \var_s) \right] \subseteq
  \left[ \exp (- 4 C \ov{S} |s|), \exp (4 C \ov{S} |s|) \right],
\end{equation}
we obtain
\begin{eqnarray}\label{eq:xxx}
    \log \int_M \exp \left( 4 (\var_s + \var_{\l,\s}) \right)
    dV_g & = & \log \int_M \exp \left( 4 \var_s \right) dV_g
    + \log \int_M \exp \left( 4 \var_{\l,\s} \right) dV_g
    \nonumber \\ & = & \log \int_M \exp \left( 4 \var_{\l,\s} \right)
    dV_g + O(\ov{S} |s|).
\end{eqnarray}
By the definition of $\var_{\l,\s}$, there holds
$$
  \int_M \exp \left( 4 \var_{\l,\s}(y) \right) d V_g(y) =
  \sum_{i=1}^k t_i \int_M \left( \frac{2 \l}{1 + \l^2 \chi_\d^2
  \left( dist(y,x_i) \right)} \right)^4 d V_g (y).
$$
We divide each of the above integrals into the metric ball $B_\d
(x_i)$ and its complement. By construction of $\chi_\d$, working in
normal coordinates centered at $x_i$, we have (for $\d$ sufficiently
small)
\begin{eqnarray*}
  & & \int_{B_\d (x_i)} \left( \frac{2 \l}{1 + \l^2 \chi_\d^2 \left(
dist(y,x_i) \right)} \right)^4 d V_g (y) = \int_{B_\d^{\R^4} (0)} (1
+ O(\d)) \left( \frac{2 \l}{1 + \l^2 |y|^2} \right)^4 d y
\\ & = & \int_{B_{\l \d}^{\R^4} (0)} (1 + O(\d)) \left( \frac{2}{1
+ |y|^2} \right)^4 d y = (1 + O(\d)) \left( \frac 83 \pi^2 + O
\left( \frac{1}{\l^4 \d^4} \right) \right).
\end{eqnarray*}
On the other hand, for $dist(y,x_i) \geq \d$ there holds $
  \left( \frac{2 \l}{1 + \l^2 \chi_\d^2 \left( dist(y,x_i)
  \right)} \right)^4 \leq \left( \frac{2 \l}{1 + \l^2
  \d^2} \right)^4.
$ Hence, from these two formulas we deduce
\begin{equation}\label{eq:intepls}
    \int_M \exp \left( 4 \var_{\l,x}(y) \right) d V_g(y) = \frac 83
    \pi^2 + O(\d) + O \left( \frac{1}{\l^4 \d^4} \right) + O
    \left( \frac{2 \l}{1 + \l^2 \d^2} \right)^4.
\end{equation}
It follows from \eqref{eq:xxx} and \eqref{eq:intepls} that
\begin{equation}\label{eq:exppls}
    \int_M \exp \left( 4 \var_{s} + 4 \var_{\l,\s} \right) dV_g =
    O(\ov{S} |s|) + O(1).
\end{equation}
This concludes the proof of \eqref{eq:estexpvar}.

\

\noindent \underline{Proof of \eqref{eq:estPvar}}. We have trivially
$$
  \langle P_g (\var_s + \var_{\l,\s}),
    (\var_s + \var_{\l,\s}) \rangle = \int_M (P_g \var_{\l,\s}, \var_{\l,\s})
dV_g + 2 \int_M (P_g \var_{s}, \var_{\l,\s}) dV_g + \int_M (P_g
\var_{s}, \var_{s}) dV_g.
$$
By Lemma \ref{l:pff} it is sufficient to estimate the last two
quantities. Since $P_g$ is negative-definite on $V$ (and since the
largest negative eigenvalue is $\l_{\ov{k}}$), we have clearly
\begin{equation}\label{eq:xxxx}
    \int_M (P_g \var_{s}, \var_{s}) dV_g \leq - |\l_{\ov{k}}|
  |s|^2 \ov{S}^2.
\end{equation}
To evaluate the second term we write $2 \int_M (P_g \var_{s},
\var_{\l,\s}) dV_g = 2 \ov{S} \sum_{i=1}^{\ov{k}} s_i \l_i \int_M
\hat{v}_i \var_{\l,\s} dV_g$. Hence it is sufficient to study each
of the terms $\int_M \hat{v}_i \var_{\l,\s} dV_g$. We claim that for
each $i \in \{1, \dots, k\}$
\begin{equation}\label{eq:intmix}
    \int_M \hat{v}_i \var_{\l,\s} dV_g = O(\d^4).
\end{equation}
In order to prove this claim, we notice first that the following
inequality holds (recall that we have chosen $\chi_\d$
non-decreasing)
$$
  \log \left( \frac{2 \l}{1 + 4 \l^2 \d^2} \right) \leq
  \var_{\l,\s} \leq \log \left( \frac{2 \l}{1 + \l^2 \chi_\d^2
  \left( d_{min}(y) \right)} \right),
$$
where $d_{min}(y) = dist(y, \{x_1\} \cup \dots \cup \{x_k\})$.
Recalling also that $\int_M \hat{v}_i dV_g = 0$, we write
$$
  \int_M \hat{v}_i(y) \var_{\l,\s}(y) dV_g(y) = \int_M \hat{v}_i
  \left( \var_{\l,\s}(y) - \log \frac{2 \l}{1 + 4 \l^2 \d^2}
  \right) dV_g.
$$
Therefore we deduce that
\begin{eqnarray*}
    \left| \int_M \hat{v}_i \var_{\l,\s} dV_g \right| &
    \leq & \|\hat{v}_i\|_{L^\infty(M)} \int_{M}
    \left( \var_{\l,\s}(y) - \log \frac{2 \l}{1 + 4 \l^2 \d^2} \right)
    dV_g(y) \\ & \leq & \|\hat{v}_i\|_{L^\infty(M)} \sum_{j=1}^k
    \int_{B_{2\d}(x_j)} \left( \log \left( \frac{2 \l}{1 + \l^2 \chi_\d^2
  \left( d_{j}(y) \right)} \right) - \log \frac{2 \l}{1 + 4 \l^2 \d^2}
  \right) dV_g(y).
\end{eqnarray*}
Working in geodesic coordinates around the point $x_j$ we get
\begin{eqnarray*}
  \int_{B_{2\d}(x_j)} \left( \var_{\l,\s}(y) - \log \frac{2 \l}{1 + 4 \l^2
  \d^2} \right)
  dV_g(y) & \leq & C \int_0^{\d} s^3 \left( \log \frac{2 \l}{1 + \l^2 s^2}
  - \log \frac{2 \l}{1 + 4 \l^2 \d^2} \right) ds \\ & + & C \int_\d^{2
  \d} s^3 \log \frac{1 + 4 \l^2 \d^2}{1 + \l^2 \chi_\d^2(s)} ds.
\end{eqnarray*}
Using elementary computations we then find
\begin{eqnarray*}\label{eq:scalar}
    \left| \int_M \hat{v}_i(y) \var_{\l,\s}(y) dV_g(y) \right| &
    \leq & C \frac{1}{\l^4} \left[ \l^4 \d^4 \log \frac{1 + 4 \l^2
    \d^2}{1 + \l^2 \d^2} + \frac{1}{8} \l^4 \d^4\right] + C \d^4 \leq C
    \d^4,
\end{eqnarray*}
which proves our claim \eqref{eq:intmix}. Notice that this
expression is independent of $\l$: this will be also used at the end
of the section. From the above formulas we obtain
$$
  \langle P_g (\var_s + \var_{\l,\s}),
    (\var_s + \var_{\l,\s}) \rangle \leq - |\l_{\ov{k}}| |s|^2 \ov{S}^2 +
    32 k \pi^2 (1 + o_\d(1)) \log \l + C_\d + O(\d^4 |s| \ov{S}),
$$
which concludes the proof of \eqref{eq:estPvar}.

\

\noindent From the three estimates \eqref{eq:estQvar},
\eqref{eq:estexpvar} and \eqref{eq:estPvar} we deduce that
\begin{eqnarray}\label{eq:IItest}
  II(\var_{\l,\s}) \leq \left( 32 k \pi^2 - 4 k_P + o_\d(1) \right)
  \log \l - |\l_{\ov{k}}| |s|^2 \ov{S}^2 + O(|s| \ov{S}) + C_\d + O(1).
\end{eqnarray}
Since $k_P > 8 k \pi^2$, choosing $\d$ sufficiently small, the
coefficient of $\log \l$ is negative. In order to show the upper
bound on $II \circ \Phi_{\ov{S},\ov{\l}}$, we fix $L > 0$. It is
easy to see that for $\ov{S}$ sufficiently large one has
$$
\left\{%
\begin{array}{ll}
    II(\var_s + 2 (1 - \var_{1,\s}) |s| + 2 \var_{1,\s} - 1) \leq -
  L & \forall \s \in M_k, \forall |s| \geq \frac
  12; \\ \\
  II(\var_s + \var_{\ov{\l},\s}) \leq - L & \forall \s \in
  M_k, \forall |s| \in \left[ \frac 14, \frac 12
  \right], \forall \l \geq 1. \\
\end{array}%
\right.
$$
After this choice of $\ov{S}$, we can also take $\ov{\l}$ so large
that
$$
    II(\var_s + \var_{\ov{\l},\s}) \leq - L, \qquad \qquad \forall
    |s| \leq \frac 14.
$$
This concludes the proof of the lemma for $k \geq 1$. In the case
$k_P < 8 \pi^2$ and $\ov{k} \geq 1$, it is sufficient to use the
estimates \eqref{eq:intps}, \eqref{eq:xx} and \eqref{eq:xxxx} to
obtain
$$
  II(\var_s) \leq - |\l_{\ov{k}}| \ov{S}^2 + O(\ov{S} |s|).
$$
The proof is thereby complete.
\end{pf}

\

\begin{pfn} {\sc of Proposition \ref{p:mkm}} The statement {\bf (a)}
follows from Lemma \ref{l:foth}. Let us prove property {\bf (b)},
starting from the case $k \geq 1$. From the expression of $e^{4
\var_{\l,\s}}$ it is easy to see that $\Psi \circ \Phi_{0,\ov{\l}}$
is homotopic to the identity on $M_k$ (to prove this, it is
sufficient to consider $\Psi \circ \Phi_{0,\l}$ for $\l$ varying
from $\ov{\l}$ to $+ \infty$). Furthermore, by continuity and by the
estimate \eqref{eq:intmix} one can check that for $|s| \leq
\frac{1}{8 \ov{S}}$ (if $\ov{S} > 1$ and if $\d$ is chosen
sufficiently small), we have
\begin{equation}\label{eq:uuu}
    \Psi (\var_s + \var_{\ov{\l},\s}) = \left(
    \hat{\Psi}(\var_s + \var_{\ov{\l},\s}), s \ov{S} + O(|s| \ov{S}
  \d^4) \right),
\end{equation}
where $\hat{\Psi}$ is defined in Lemma \ref{l:defmap}, and therefore
$\Psi \circ \Phi_{\ov{S},\ov{\l}}$ is homotopic (in $A_{k,\ov{k}}$)
to the identity on $M_k \times B_{\frac{1}{8 \ov{S}}}^{\ov{k}}
\subseteq A_{k,\ov{k}}$.


On the other hand, by \eqref{eq:uuu}, for $|s| \geq \frac{1}{8
\ov{S}}$, the $\ov{k}$-vector $s \ov{S} + O(|s| \ov{S} \d^4)$ is
almost parallel to $s$ (and non-zero), and therefore on this set
$\Psi \circ \Phi_{\ov{S},\ov{\l}}$ can be easily contracted to the
boundary of $B_1^{\ov{k}}$ (recall the definition of
$A_{k,\ov{k}})$, as for the identity map. This concludes the proof
in the case $k \geq 1$. The proof for $k_P < 8 \pi^2$ and under the
assumption \eqref{eq:kp4} is analogous.
\end{pfn}

\section{Proof of Theorem \ref{th:ex}}\label{s:proof}

In this section we prove Theorem \ref{th:ex} employing a min-max
scheme based on the construction of the above set $A_{k,\ov{k}}$,
see Lemma \ref{l:min-max}. As anticipated in the introduction, we
then define a modified functional $II_{\rho}$ for which we can prove
existence of solutions in a dense set of the values of $\rho$.
Following an idea of Struwe, this is done proving the a.e.
differentiability of the map $\rho \mapsto \ov{\Pi}_\rho$, where
$\ov{\Pi}_\rho$ is the min-max value for the functional $II_\rho$.

\

\noindent We now introduce the scheme which provides existence of
solutions for \eqref{eq:Qc}, beginning with the case $k \geq 1$. Let
$\widehat{A_{k,\ov{k}}}$ denote the (contractible) cone over
$A_{k,\ov{k}}$, which can be represented as $\widehat{A_{k,\ov{k}}}
= A_{k,\ov{k}} \times [0,1]$ with $A_{k,\ov{k}} \times \{0\}$
collapsed to a single point. Let first $L$ be so large that
Proposition \ref{p:map} applies with $\frac L4$, and then let
$\ov{S}, \ov{\l}$ be so large (and $\d$ so small) that Proposition
\ref{p:mkm} applies for this value of $L$. Fixing these numbers
$\ov{S}$ and $\ov{\l}$, we define the following class
\begin{equation}\label{eq:PiPi}
    \Pi_{\ov{S},\ov{\l}} = \left\{ \pi : \widehat{A_{k,\ov{k}}} \to H^2(M)
 \; : \; \pi \hbox{ is continuous and } \pi(\cdot \times \{1\})
 = \Phi_{\ov{S},\ov{\l}}(\cdot) \hbox{ on } A_{k,\ov{k}} \right\}.
\end{equation}
In the case $k_P < 8 \pi^2$ and $\ov{k} \geq 1$ we simply use the
closed unit $\ov{k}$-dimensional ball $\ov{B}^{\ov{k}}_1$ and we
set (still for large values of $L$)
$$
 \Pi_{\ov{S}} = \left\{ \pi : \ov{B}^{\ov{k}}_1 \to H^2(M)
 \; : \; \pi \hbox{ is continuous and } \pi(\cdot)
 = \Phi_{\ov{S}}(\cdot) \hbox{ on } S^{\ov{k}-1} \right\}.
$$
Then we have the following properties.

\begin{lem}\label{l:min-max}
The set $\Pi_{\ov{S},\ov{\l}}$ (resp. $\Pi_{\ov{S}}$) is non-empty
and moreover, letting
$$
  \ov{\Pi}_{\ov{S},\ov{\l}} = \inf_{\pi \in \Pi_{\ov{S},\ov{\l}}}
  \; \sup_{m \in \widehat{A_{k,\ov{k}}}} II(\pi(m)), \qquad
  \hbox{ there holds } \qquad \ov{\Pi}_{\ov{S},\ov{\l}} > - \frac
  L2,
$$
$$
  \left(resp. \;\; \ov{\Pi}_{\ov{S}} = \inf_{\pi \in \Pi_{\ov{S}}}
  \; \sup_{m \in \ov{B}^{\ov{k}}_1} II(\pi(m)), \qquad
  \hbox{ there holds } \qquad \ov{\Pi}_{\ov{S}} > -
  \frac L2 \right).
$$
\end{lem}

\begin{pf}
To prove that $\Pi_{\ov{S},\ov{\l}} \neq \emptyset$, we just
notice that the following map
\begin{equation}\label{eq:ovPi}
  \ov{\pi}(z,t) = t \Phi_{\ov{S},\ov{\l}} (z), \qquad \qquad (z,t)
  \in \widehat{A_{k,\ov{k}}},
\end{equation}
belongs to $\Pi_{\ov{S},\ov{\l}}$. Assuming by contradiction that
$\ov{\Pi}_{\ov{S},\ov{\l}} \leq  - \frac L2$, there would exist a
map $\pi \in \Pi_{\ov{S},\ov{\l}}$ with $\sup_{m \in
\widehat{A_{k,\ov{k}}}} II(\pi(m)) \leq - \frac 38 L$. Then, since
Proposition \ref{p:map} applies with $\frac L4$, writing $m = (z,
t)$, with $z \in A_{k,\ov{k}}$, the map
$$
  t \mapsto \Psi \circ \pi(\cdot,t)
$$
would be an homotopy in $A_{k,\ov{k}}$ between $\Psi \circ
\Phi_{\ov{S},\ov{\l}}$ and a constant map. But this is impossible
since $A_{k,\ov{k}}$ is non-contractible (see Corollary \ref{c:akk})
and since $\Psi \circ \Phi_{\ov{S},\ov{\l}}$ is homotopic to the
identity on $A_{k,\ov{k}}$, by Proposition \ref{p:mkm}. Therefore we
deduce $\ov{\Pi}_{\ov{S},\ov{\l}} > - \frac L2$.

In the case $k_P < 8 \pi^2$ and $\ov{k} \geq 1$ it is sufficient to
take $\ov{\pi}(z,t) = t \Phi_{\ov{S}} (z)$ and to proceed in the
same way.
\end{pf}

\

\noindent Next we introduce a variant of the above min-max scheme,
following \cite{str} and \cite{djlw}. When $k_P < 8 \pi^2$, we
define for convenience $A_{k,\ov{k}} = S^{\ov{k}}$,
$\widehat{A_{k,\ov{k}}} = \ov{B}^{\ov{k}}_1$, $\Phi_{\ov{S},\ov{\l}}
= \Phi_{\ov{S}}$, etc. For $\rho$ in a small neighborhood of $1$, $[
1 - \rho_0, 1 + \rho_0]$, we define the modified functional $II_\rho
: H^2(M) \to \R$ as
\begin{equation}\label{eq:IIrho}
    II_\rho(u) = \langle P_g u, u \rangle + 4 \rho \int_M Q_g u -
    4 \rho k_P \log \int_M e^{4 u} dV_g.
\end{equation}
Following the estimates of the previous sections, one easily checks
that the above min-max scheme applies uniformly for $\rho \in [ 1 -
\rho_0, 1 + \rho_0 ]$ and for $\ov{S}$, $\ov{\l}$ sufficiently
large. More precisely, given any large number $L > 0$, there exist
$\rho_0$ sufficiently small and $\ov{S}$, $\ov{\l}$ sufficiently
large so that
\begin{equation}\label{eq:min-maxrho}
   \sup_{\pi \in {\Pi}_{\ov{S},\ov{\l}}} \sup_{m \in \partial
   \widehat{A_{k,\ov{k}}}} II_{\rho}(\pi(m)) < - 2 L; \quad
   \ov{\Pi}_{\rho} := \inf_{\pi \in \Pi_{\ov{S},\ov{\l}}}
  \; \sup_{m \in \widehat{A_{k,\ov{k}}}} II(\pi(m)) > -
  \frac{L}{2}; \qquad  \rho \in [1 - \rho_0, 1 + \rho_0],
\end{equation}
where $\Pi_{\ov{S},\ov{\l}}$ is defined in \eqref{eq:PiPi}.
Moreover, using for example the test map \eqref{eq:ovPi}, one
shows that for $\rho_0$ sufficiently small there exists a large
constant $\ov{L}$ such that
\begin{equation}\label{eq:ovlovl}
  \ov{\Pi}_{\rho} \leq \ov{L} \qquad \qquad \hbox{ for every }
  \rho \in [1 - \rho_0, 1 + \rho_0].
\end{equation}

\

\noindent We have the following result, regarding the dependence in
$\rho$ of the min-max value $\ov{\Pi}_\rho$. A similar statement has
been proved in \cite{djlw}, but here we allow the presence of
negative eigenvalues for the elliptic operator, so the proof is more
involved. Since this is rather technical, we give it in the
appendix.

\begin{lem}\label{l:arho}
Let $\ov{S}$, $\ov{\l}$ be so large and $\rho_0$ be so small that
\eqref{eq:min-maxrho} holds. Then, taking $\rho_0$ possibly smaller,
there exists a fixed constant $C$ (depending only on $M$ and
$\rho_0$) such that the function
$$
  \rho \mapsto \frac{\ov{\Pi}_\rho}{\rho} - C \rho \qquad \hbox{ is
  non-increasing in } [1 - \rho_0, 1 + \rho_0].
$$
\end{lem}

\

\noindent From Lemma \ref{l:arho} we deduce that the function $\rho
\mapsto \frac{\ov{\Pi}_\rho}{\rho}$ is differentiable a.e., and we
obtain the following corollary.

\begin{cor}\label{c:c}
Let $\ov{S}$, $\ov{\l}$ and $\rho_0$ be as in Lemma \ref{l:arho},
and let $\L \subset [1 - \rho_0, 1 + \rho_0]$ be the (dense) set
of $\rho$ for which the function $\frac{\ov{\Pi}_\rho}{\rho}$ is
differentiable. Then for $\rho \in \L$ the functional $II_\rho$
possesses a bounded Palais-Smale sequence $(u_l)_l$ at level
$\ov{\Pi}_{\rho}$.
\end{cor}

\begin{pf}
The existence of a Palais-Smale sequence $(u_l)_l$ follows from
Lemma \ref{l:min-max}, and the boundedness is proved exactly as in
\cite{djlw}, Lemma 3.2.
\end{pf}

\begin{rem}\label{r:k<8pi2}
When $k_P < 8 \pi^2$ one can use a direct approach to prove
boundedness of Palais-Smale sequences (satisfying
\eqref{eq:noul}). We test the relation $II'(u_l) \to 0$ (in
$H^{-2}(M)$) on $\hat{u}_l$ and $\tilde{u}_l$, where $\hat{u}_l$
is the component of $u_l$ in $V$ and $\tilde{u}_l$ is the
component perpendicular to $V$.

Testing on $\hat{u}_l$ we obtain
\begin{equation}\label{eq:uuuu}
    \langle P_g \hat{u}_l, \hat{u}_l \rangle + 4 \int_M Q_g
    \hat{u}_l dV_g - 4 k_P \int_M
  e^{4 u_l} \hat{u}_l dV_g = o_l(1) \|\hat{u}_l\|_{L^\infty(M)}.
\end{equation}
Since $\|e^{4 u_l}\|_{L^1(M)} = 1$ by \eqref{eq:noul} and since on
$V$ the $L^\infty$-norm is equivalent to the $H^2$-norm, the last
formula implies $- \langle P_g \hat{u}_l, \hat{u}_l \rangle = O(1)
\|\hat{u}_l\|_{H^2(M)}$. Therefore, being $P_g$ negative-definite on
$V$, we get uniform bounds on $\|\hat{u}_l\|$.

On the other hand, testing the equation on $\tilde{u}_l$ we find
$$
   2 \langle P_g \tilde{u}_l, \tilde{u}_l \rangle - 4 k_P \int_M
   e^{4 u_l} (\tilde{u}_l - \ov{u}_l) dV_g = O(\| \tilde{u}_l - \ov{u}_l
   \|_{H^2(M)}).
$$
This implies, for any $\a > 1$ (using \eqref{eq:adaneg} and
\eqref{eq:uuuu})
\begin{eqnarray*}
    2 \langle P_g \tilde{u}_l, \tilde{u}_l \rangle & \leq & C
    e^{4 \ov{u}_l} \int_M e^{4 (u_l - \ov{u}_l)} (\tilde{u}_l -
    \ov{u}_l) dV_g + O(\| \tilde{u}_l - \ov{u}_l \|_{H^2(M)}) \\
    & \leq & C_\a
    e^{4 \ov{u}_l} \int_M e^{4 \a (u_l - \ov{u}_l)} dV_g +
    O(\| \tilde{u}_l - \ov{u}_l \|_{H^2(M)}) \\ & \leq &
    C_\a e^{4 \ov{u}_l} e^{\a^2 \frac{\langle P_g \tilde{u}_l,
    \tilde{u}_l \rangle}{8 \pi^2}} +
    O(\| \tilde{u}_l - \ov{u}_l \|_{H^2(M)}).
\end{eqnarray*}
Moreover, since we are assuming $II(u_l) \to c \in \R$, for any
small $\e$ we get
\begin{eqnarray*}
    C \geq II(\tilde{u}_l) = \langle P_g \tilde{u}_l,
    \tilde{u}_l \rangle + 4 \int_M Q_g \tilde{u}_l =
    (1 + O(\e)) \langle P_g \tilde{u}_l,
    \tilde{u}_l \rangle + 4 k_P \ov{u}_l + C_\e,
\end{eqnarray*}
provided $l$ is sufficiently large. Hence from the last two
formulas we deduce
$$
  \langle P_g \tilde{u}_l, \tilde{u}_l \rangle \leq C_{\a,\e}
  e^{\langle P_g \tilde{u}_l, \tilde{u}_l \rangle \left(
  \frac{\a^2}{8 \pi^2} - \frac{1+O(\e)}{k_P} \right)} +
  O(\| \tilde{u}_l - \ov{u}_l \|_{H^2(M)}).
$$
Now, choosing $\a$ close to $1$ and $\e$ so small that the
exponential factor has a negative coefficient (this is always
possible since $k_P < 8 \pi^2$), we obtain a uniform bound for
$\|\tilde{u}_l - \ov{u}_l\|$. The bound on $\ov{u}_l$ now follows
easily from \eqref{eq:noul}.
\end{rem}

\

\noindent  Now the proof of Theorem \ref{th:ex} is an easy
consequence of the following proposition and of Theorem \ref{th:bd}.

\begin{pro}\label{p:weakstr}
Suppose $(u_l)_l \subseteq H^2(M)$ is a sequence for which (as $l
\to + \infty$)
$$
  II_\rho(u_l) \to c \in \R; \qquad \qquad II'_\rho(u_l) \to 0;
  \qquad \qquad \|u_l\|_{H^2(M)} \leq C,
$$
where $C$ is independent of $l$. Then $(u_l)_l$ has a weak limit
$u_0$ which satisfies \eqref{eq:mod}.
\end{pro}

\begin{pf} The existence of a weak limit $u_0 \in H^2(M)$ follows from
Corollary \ref{c:c}. Let us show that $u_0$ satisfies $II'_\rho(u_0)
= 0$. For any function $v \in H^2(M)$ there holds
\begin{eqnarray*}
    II'_\rho(u_0)[v] & = & II'_\rho(u_l)[v] + 2 \langle P_g v, (u_0 - u_l)
    \rangle + 4 \rho k_P \left( \frac{\int_M e^{4 u_l} v dV_g}{\int_M e^{4 u_l} dV_g}
    - \frac{\int_M e^{4 u_0} v dV_g}{\int_M e^{4 u_0} dV_g} \right).
\end{eqnarray*}
Since the first two terms in the right-hand side tend to zero by our
assumptions, it is sufficient to check that $\int_M e^{4 u_l} v dV_g
= \int_M e^{4 u_0} v dV_g + o(1) \|v\|_{H^2(M)}$ (to deal with the
denominators just take $v \equiv 1$). In order to do this, consider
$p, p', p''
> 1$ satisfying $\frac{1}{p} + \frac{1}{p'} + \frac{1}{p''} = 1$.
Using Lagrange's formula we obtain, for some function $\th_l$ with
range in $[0,1]$, $e^{4 u_l} - e^{4 u_0} = e^{4\th_l u_l + 4(1 -
\th_l)u_0} (u_l - u_0)$ a.e. in $x$. Then from some elementary
inequalities we find
\begin{eqnarray*}
    \int_M \left( e^{4 u_l} - e^{4 u_0} \right) v \, dV_g & \leq &
    C \int_M \left( e^{4 u_l} + e^{4 u_0} \right) |u_l -
    u_0| |v| dV_g \\ & \leq & C \left[ \|e^{4 u_l}\|_{L^p(M)} +
    \|e^{4 u_0}\|_{L^p(M)} \right] \|u_l - u_0 \|_{L^{p'}(M)}
    \|v\|_{L^{p''}(M)} \\ & \leq & o(1) \|v\|_{L^{p''}(M)} = o(1) \|v\|_{H^2(M)},
\end{eqnarray*}
by \eqref{eq:adaneg}, the boundedness of $(u_l)_l$ and the fact that
$u_l \rightharpoonup u_0$ weakly in $H^2(M)$.
\end{pf}

\section{Appendix}

In this section we collect the most technical proofs of the paper,
namely those of Lemmas \ref{l:Pjl}, \ref{l:pff} and \ref{l:arho}.

\

\begin{pfnb} {\bf of Lemma \ref{l:Pjl}} \quad
By Corollary \ref{cor:Mksmo}, we know that $M_j\left(
\frac{\e}{4}\right)$ is a smooth finite-dimensional manifold.
Therefore, if $\hat{\e}$ is sufficiently small, there exists a
continuous projection $P_j$ from $\hat{M}_{k,j}^{\hat{\e},\e}$ onto
$M_j\left( \frac{\e}{2} \right)$ (whose closure lies in $M_j\left(
\frac{\e}{4}\right)$). Since we are regarding $M_k$ as a subset of
$C^1(M)^*$, a Banach space, we cannot in general project elements in
a neighborhood of $M_j\left( \frac{\e}{2} \right)$ onto their
closest point in $M_j\left( \frac{\e}{2} \right)$ (this might not be
unique). Nevertheless, using the Implicit Function Theorem and a
partition of unity it is possible to define the projection in such a
way that
\begin{equation}\label{eq:projp}
    dist(\s,P_j(\s)) \leq C_{k,\e} dist(\s, M_j\left( \e \right)),
    \qquad \quad \s \in  \hat{M}_{k,j}^{\hat{\e},\e},
\end{equation}
where $C_{k,\e}$ is a constant depending only on $k$ and $\e$ (we
are taking $1 \leq j \leq k-1$).

To fix some notation, we use the following convention
$$
  \s = \sum_{i=1}^k t_i \d_{x_i}; \qquad \qquad P_j (\s) =
  \sum_{i=1}^j s_i \d_{y_i}.
$$
By Lemma \ref{l:Pj-1}, since we are assuming that $P_j(\s)$ belongs
to $M_j(\frac \e 2)$, we have the following estimates
$$
  s_i \geq \frac \e4, \qquad \qquad dist(y_i, y_l) \geq \frac{\e}{4};
  \qquad \qquad \forall \; \, i, l = 1, \dots, s, i \neq l.
$$
Moreover the points $y_i$ and the numbers $s_i$ depend continuously
on $\s$.

We define first an auxiliary map $\tilde{T}^t_j$, $\tilde{T}^t_j(\s)
= \sum \tilde{t}_i \d_{\tilde{x}_i}$, which misses the normalization
condition $\sum_{i=1}^k \tilde{t}_i = 1$, but only up to a small
error. This map will then be corrected to the real $T^t_j$. The idea
to construct $\tilde{T}^t_j$ is the following. If a point $x_i$ is
far from each $y_l$, we keep this point fixed and let its
coefficient vanish to zero as $t$ varies from 0 to 1. On the other
hand, if $x_i$ is close to some of the $y_l$'s, then we translate it
to a {\em weighted convex combination}  of the points $x_i$ which
are close to the same $y_l$.

To make this construction rigorous (and the map $\tilde{T}^t_j$
continuous), we consider a small number $\eta \ll \e$ (this will be
chosen later of order $C_{k,\e} \sqrt{\hat{\e}}$) (where $C_{k,\e}$
depends only on $k$ and $\e$), and define a smooth cutoff function
$\rho_\eta$ satisfying the following properties
\begin{equation}\label{eq:variPj}
    \left\{%
\begin{array}{ll}
    \rho_\eta(t) = 1, & \hbox{ for } t \leq \frac{\eta}{16};  \\
    \rho_\eta(t) = 0, & \hbox{ for } t \geq \frac{\eta}{8}; \\
    \rho_\eta(t) \in [0,1], & \hbox{ for every } t \geq 0. \\
\end{array}%
\right.
\end{equation}
Then we set
\begin{equation}\label{eq:rhoi}
    \rho_{l,\eta}(x) = \rho_\eta(dist(x,y_l)); \qquad \qquad \hbox{ for
    } l = 1, \dots, j.
\end{equation}
We define also the following quantities
$$
  \mathcal{T}_l(\s) = \sum_{x_i \in B_{\frac \eta 8}(y_l)}
  \rho_{l,\eta}(x_i) t_i; \qquad \mathcal{X}_l(\s) = \frac{1}{\mathcal{T}_l(\s)}
  \sum_{x_i \in B_{\frac \eta 8}(y_l)} \rho_{l,\eta}(x_i) t_i x_i,
  \qquad \quad l = 1, \dots, j.
$$
We notice that, if $\eta$ is chosen sufficiently small, the weighted
convex combination $\mathcal{X}_l(\s)$ is well-defined, see the
notation in Section \ref{s:pr}. We also set
$$
  z_i(\s) = \frac{8}{\eta} dist(x_i, y_l) - 1, \qquad \qquad \hbox{ for }
  x_i \in B_{\frac \eta 4}(y_l).
$$
Since for all $i \neq l$ there holds $dist(y_i, y_l) \geq \frac \e
4$ and since $\eta << \e$, then for every $i$ there exists at most
one point $y_l$ such that $x_i \in B_{\frac \eta 4}(y_l)$, and hence
the number $z_i(\s)$ is well defined. Now we construct the map
$\tilde{T}^t_j(\s)$ as follows
$$
  \tilde{T}^t_j(\s) = \sum_{i=1}^k \tilde{t}_i(\s,t)
  \d_{\tilde{x}_i(\s,t)},
$$
where the numbers $\tilde{t}_i(\s,t)$ and the points
$\tilde{x}_i(\s,t)$ are given by
\begin{eqnarray*}
  \tilde{t}_i(\s,t) = (1-t) t_i; \qquad \tilde{x}_i(\s,t) = x_i &
  \hbox{ if } & x_i \in M \setminus \cup_l B_{\frac \eta 4}(y_l); \\
  \tilde{t}_i(\s,t) = (1-t) t_i; \qquad \tilde{x}_i(\s,t) = (1-t)
  x_i + t[z_i(\s) x_i + (1 - z_i(\s)) \mathcal{X}_l(\s)] &
  \hbox{ if } & x_i \in B_{\frac \eta 4}(y_l) \setminus B_{\frac \eta 8}(y_l); \\
  \tilde{t}_i(\s,t) = ((1-t) + t \rho_{l,\eta}(x_i)) t_i; \qquad
  \tilde{x}_i(\s,t) = (1 - t) x_i + t \mathcal{X}_l(\s) &
  \hbox{ if } & x_i \in B_{\frac \eta 8}(y_l).
\end{eqnarray*}
As already mentioned, the numbers $\tilde{t}_i(\s,t)$ will in
general miss the condition $\sum_i \tilde{t}_i(\s,t) = 1$. The next
step consists in estimating this sum and correct the map
$\tilde{T}^t_j(\s)$ in order to match this condition. For this
purpose it is convenient to define
$$
  \tilde{\mathcal{T}}_l(\s,t) = \sum_{x_i \in
  B_{\frac{\eta}{8}}(y_l)} \tilde{t}_i(\s,t); \qquad \qquad
  \tilde{\mathcal{T}}(\s,t) = 1 - \sum_{l=1}^j \tilde{\mathcal{T}}_l(\s,t).
$$
Now we finally set
\begin{equation}\label{eq:T^t_j}
    T^t_j(\s) = \frac{1}{(1-t)\tilde{\mathcal{T}}(\s,0) +
    \sum_{l=1}^j \tilde{\mathcal{T}}_l(\s,t)} \sum_{i=1}^k \tilde{t}_i(\s,t)
  \d_{\tilde{x}_i(\s,t)}.
\end{equation}
We notice that the sum of all the coefficients is $1$, and that the
map is well defined and continuous in both $t$ and $\s$. We also
notice that the properties {\bf (i)}, {\bf (iv)} and {\bf (v)} are
satisfied, while {\bf (ii)} follows from {\bf (iii)}. Therefore  it
only remains to prove {\bf (iii)}. First of all we give an estimate
on the quantities $\tilde{\mathcal{T}}_l(\s,t)$ and
$\tilde{\mathcal{T}}(\s,t)$.

We recall that we have taken $\s \in \hat{M}^{\hat{\e},\e}_{k,j}$,
and hence by \eqref{eq:projp} for any function $f \in C^1(M)$ with
$\|f\|_{C^1(M)} \leq 1$ one has $\left| (\s - P_j(\s), f) \right|
\leq C_{k,\e} \hat{\e}$. We now choose a function $f$ satisfying the
following properties
$$
  f(x) = \left\{%
\begin{array}{ll}
    \frac 12 & \hbox{ for } x \in \cup_{l=1}^j B_{\frac{\eta}{48}}(y_l); \\
    \frac 12 + \frac{\eta}{32} & \hbox{ for } x \in M \setminus
    \cup_{l=1}^j B_{\frac{\eta}{16}}(y_l);  \\
    \|f\|_{C^1(M)} \leq 1. & \\
\end{array}%
\right.
$$
For this function we have $(P_j(\s), f) = \sum_{i=1}^j s_i f(y_i) =
\frac 12$ and moreover
\begin{eqnarray*}
  (\s, f) & = & \sum_{x_i \in \cup_{l=1}^j B_{\frac{\eta}{16}}(y_l)} t_i
  f(x_i) + \sum_{x_i \in M \setminus \cup_{l=1}^j B_{\frac{\eta}{16}}(y_l)} t_i
  f(x_i) \\ & \geq & \frac{1}{2} \sum_{x_i \in \cup_{l=1}^j B_{\frac{\eta}{16}}(y_l)}
  t_i + \left( \frac 12 + \frac{\eta}{32} \right) \sum_{x_i \in M \setminus
  \cup_{l=1}^j B_{\frac{\eta}{16}}(y_l)} t_i.
\end{eqnarray*}
Therefore we deduce the following inequality
\begin{equation*}
    \frac{\eta}{32} \sum_{x_i \in M \setminus \cup_{l=1}^j B_{\frac{\eta}{16}}(y_l)}
    t_i \leq (\s, f) - (P_j(\s), f) \leq C_{k,\e} \hat{\e}.
\end{equation*}
This estimate implies
$$
   \tilde{\mathcal{T}}(\s,0) = \sum_{x_i \in M \setminus \cup_{l=1}^j
   B_{\frac{\eta}{16}}(y_l)} t_i \leq 32 \frac{C_{k,\e} \hat{\e}}{\eta},
$$
and also (since $\rho_{l,\eta} \equiv 1$ in
$B_{\frac{\eta}{16}}(y_l)$)
\begin{eqnarray*}
    \tilde{\mathcal{T}}_l(\s,t) & = & \sum_{x_i \in
    B_{\frac{\eta}{8}}(y_l) \setminus B_{\frac{\eta}{16}}(y_l)} \left( (1-t) + t
    \rho_{l,\eta}(x_i)
    \right) t_i + \sum_{x_i \in B_{\frac{\eta}{16}}(y_l)} \left(
    (1-t) + t \rho_{l,\eta}(x_i) \right) t_i \\ & = &
    \tilde{\mathcal{A}}_l(\s,t) + \sum_{x_i \in B_{\frac{\eta}{16}}(y_l)}
    t_i, \qquad \qquad \hbox{ where } \quad \sum_{l=1}^j
    |\tilde{\mathcal{A}}_l(\s,t)| \leq 32 \frac{C_{k,\e} \hat{\e}}{\eta}.
\end{eqnarray*}
Hence, since $\sum_{l=1}^j \tilde{\mathcal{T}}_l(\s,0) +
\tilde{\mathcal{T}}(\s,0) = 1$, there holds
$$
  1 = \sum_{x_i \in \cup_{l=1}^j B_{\frac{\eta}{16}}(y_l)} t_i + \sum_{l=1}^j
  \tilde{\mathcal{A}}_l(\s,0) + \sum_{x_i \in M \setminus \cup_{l=1}^j
  B_{\frac{\eta}{8}}(y_l)} t_i,
$$
from which we deduce
\begin{eqnarray*}
    \left| \sum_{l=1}^j \tilde{\mathcal{T}}_l(\s,t) + (1-t) \tilde{\mathcal{T}}(0)
    - 1 \right| & = & \left| \left( \sum_{l=1}^j \left(\tilde{\mathcal{A}}_l(\s,t) -
    \tilde{\mathcal{A}}_l(\s,0)\right) \right) + (1-t) \sum_{x_i \in M \setminus \cup_{l=1}^j
  B_{\frac{\eta}{8}}(y_l)} t_i \right| \\ & \leq & 64
  \frac{C_{k,\e} \hat{\e}}{\eta} + 32  \frac{C_{k,\e} \hat{\e}}{\eta} = 96
  \frac{C_{k,\e} \hat{\e}}{\eta}.
\end{eqnarray*}
As a consequence, using a Taylor expansion (recall that we are
choosing $\frac{C_{k,\e} \hat{\e}}{\eta} \ll 1$), we find that the
coefficient added in the definition of $T^t_j$, see
\eqref{eq:T^t_j}, can be estimated by
$$
  \left| \frac{1}{\sum_{l=1}^j \tilde{\mathcal{T}}_l(\s,t) + (1-t)
  \tilde{\mathcal{T}}(0)} - 1 \right| \leq 100 \frac{C_{k,\e} \hat{\e}}{\eta}.
$$
To control the metric distance in {\bf (iii)}, we use the last
formula to get, for an arbitrary function $f \in C^1(M)$ with
$\|f\|_{C^1(M)} \leq 1$
\begin{eqnarray}\label{eq:tttilde}
    \left| (\s, f) - (T^t_j(\s), f) \right| & \leq & \left| (\s, f)
    - (\tilde{T}^t_j(\s), f) \right| + \left| (\tilde{T}^t_j(\s), f)
    - (T^t_j(\s), f) \right| \nonumber \\ & \leq & \left| (\s, f)
    - (\tilde{T}^t_j(\s), f) \right| + 100 \frac{C_{k,\e} \hat{\e}}{\eta}.
\end{eqnarray}
Hence it is sufficient to estimate the distance between $\s$ and
$\tilde{T}^t_j(\s)$. We can write
\begin{eqnarray*}
  \left| (\s, f) - (\tilde{T}^t_j(\s), f) \right| & \leq & \sum_{x_i \in M
  \setminus \cup_{l=1}^j B_{\frac{\eta}{4}}(y_l)} t_i + \sum_{x_i \in \cup_{l=1}^j
  (B_{\frac{\eta}{4}}(y_l) \setminus B_{\frac{\eta}{16}}(y_l))} \left|
   t_i f(x_i) - \tilde{t}_i(\s,t) f(\tilde{x}_i(\s,t)) \right| \\
   & + & \sum_{x_i \in \cup_{l=1}^j B_{\frac{\eta}{16}}(y_l)} t_i
   dist(x_i, \tilde{x}_i(\s,t)).
\end{eqnarray*}
Since $\left| t_i f(x_i) - \tilde{t}_i(\s,t) f(\tilde{x}_i(\s,t))
\right| \leq |t_i - \tilde{t}_i(\s,t)| + \tilde{t}_i(\s,t) dist(x_i,
\tilde{x}_i(\s,t)) \leq 2 t_i$ (for $\eta$ small), we obtain
\begin{eqnarray*}
  \left| (\s, f) - (\tilde{T}^t_j(\s), f) \right| & \leq & 2 \sum_{x_i
  \in M \setminus \cup_{l=1}^j B_{\frac{\eta}{16}}(y_l)} t_i + \sum_{l=1}^j \sum_{x_i \in
  B_{\frac{\eta}{16}}(y_l)} t_i dist(x_i, \mathcal{X}_l(\s)) \\ & \leq &
  64 \frac{C_{k,\e} \hat{\e}}{\eta} + \sum_{l=1}^j \sum_{x_i \in
  B_{\frac{\eta}{16}}(y_l)} t_i dist(x_i, \mathcal{X}_l(\s)).
\end{eqnarray*}
In order to estimate the last term, we notice that each point $x_i$
in the homotopy is shifted at most of $\frac{\eta}{2}$, see the
comments at the beginning of Section \ref{s:pr}. Therefore from
\eqref{eq:tttilde} and the last expression we get
$$
  \left| (\s, f) - (T^t_j(\s), f) \right| \leq 170
  \frac{C_{k,\e} \hat{\e}}{\eta} + \frac{\eta}{2}.
$$
Therefore, choosing $\eta = C_{k,\e} \sqrt{\hat{\e}}$, we obtain the
desired conclusion.
\end{pfnb}

\

\begin{pfnb} {\bf of Lemma \ref{l:pff}} \quad For simplicity, see Section
\ref{s:test}, we adopt again the notation $d_i = d_i(y) = dist(y,
x_i)$, $i = 1, \dots, k$, and we consider these as functions of $y$,
for $(x_i)_i$ fixed. By \eqref{eq:pls} with some straightforward
computations we find
\begin{equation}\label{eq:njvar}
    \n \var_{\l,\s} = - \l^2 (2 \l)^4 \frac{\sum_{i=1}^k t_i \n
  (\chi_\d^2(d_i)) (1 + \l^2 \chi_\d^2(d_i))^{-5}}{\sum_{s=1}^k t_s
  \left( \frac{2 \l}{1 + \l^2 \chi_\d^2(d_s)} \right)^4},
\end{equation}
and
\begin{eqnarray}\label{eq:Dvar}
    \D \var_{\l,\s} & = & \l^2 (2 \l)^4 \frac{\sum_{i=1}^k t_i
    (1 + \l^2 \chi_\d^2(d_i))^{-6} \left[ 5 \l^2 |\n (\chi_\d^2(d_i))|^2
    - \D (\chi_\d^2(d_i)) (1 + \l^2 \chi_\d^2(d_i)) \right]}{\sum_{s=1}^k
    t_s \left( \frac{2 \l}{1 + \l^2 \chi_\d^2(d_s)} \right)^4} \nonumber
    \\ & - & 4 \l^4 (2 \l)^8 \frac{\sum_{i,s=1}^k t_i t_s (1 + \l^2
    \chi_\d^2(d_i))^{-5} (1 + \l^2 \chi_\d^2(d_s))^{-5} \n (\chi_\d^2 (d_i))
    \cdot \n (\chi_\d^2 (d_s))}{\left[ \sum_{r=1}^k t_r \left( \frac{2 \l}{1
    + \l^2 \chi_\d^2(d_r)} \right)^4 \right]^2}.
\end{eqnarray}
We begin by estimating $\int_M (\D \var_{\l,\s})^2 dV_g$. This is
the most involved part of the proof, and the result is given in
formula \eqref{eq:La2m} below. We notice first that the following
pointwise estimate holds true, as one can easily check using
\eqref{eq:Dvar}
$$
  |\D \var_{\l,\s}| \leq \frac{C}{\l^2}.
$$
For a large but fixed constant $\Theta > 0$ (and for $\l \to +
\infty$), the volume of a ball in $M$ of radius $\frac{\Theta}{\l}$
is bounded by $C \frac{\Theta^4}{\l^4}$. From this we deduce that
\begin{equation}\label{eq:intL2Bi}
    \int_{\cup_{i=1}^k B_{\frac{\Theta}{\l}}(x_i)} (\D \var_{\l,\s})^2 dV_g \leq C
  \Theta^4.
\end{equation}
Therefore we just need to estimate the integral on the complement of
the union of these balls, which we denote by
\begin{equation}\label{eq:MsT}
    M_{\s,\Theta} = M \setminus \cup_{i=1}^k
    B_{\frac{\Theta}{\l}}(x_i).
\end{equation}
In this set, since we are taking $\Theta$ large, the ratio between
$1 + \l^2 d_i^2$ and $\l^2 d_i^2$ is very close to 1, and hence we
obtain the following estimates
\begin{equation}\label{eq:appth}
    (1 + \l^2 \chi_\d^2(d_i)) = (1 + o_{\d,\Theta}(1)) \l^2
  \chi_\d^2(d_i) \qquad \qquad \hbox{ in } M_{\s,\Theta};
\end{equation}
\begin{equation}\label{eq:appth2}
    5 \l^2 |\n (\chi_\d^2(d_i))|^2 - \D (\chi_\d^2(d_i)) (1 + \l^2
  \chi_\d^2(d_i)) = 12 (1 + o_{\d,\Theta}(1)) \l^2 \tilde{\chi}_\d^2(d_i)
  \qquad \qquad \hbox{ in } M_{\s,\Theta},
\end{equation}
where $o_{\d,\Theta}(1)$ tends to zero as $\d$ tends to zero and
$\Theta$ tends to infinity, and where $\tilde{\chi}_\d$ is a new
cutoff function (which depends on $\chi_\d$) satisfying
\begin{equation}\label{eq:tchid}
    \left\{%
\begin{array}{ll}
    \tilde{\chi}_\d(t) =  t, & \hbox{ for } t \in [0,\d]; \\
    \tilde{\chi}_\d(t) = 0, & \hbox{ for } t \geq 2 \d; \\
    \tilde{\chi}_\d(t) \in [0, 2 \d], & \hbox{ for } t \in [\d, 2 \d]. \\
\end{array}%
\right.
\end{equation}
Using \eqref{eq:Dvar}, \eqref{eq:appth} and \eqref{eq:appth2} one
finds that the following estimate holds
\begin{eqnarray}\label{eq:bettest}
    \D \var_{\l,\d} & = & 12 (1 + o_{\d,\Theta}(1))
    \frac{\sum_{i=1}^k t_i
    \frac{\tilde{\chi}_\d^2(d_i)}{\chi_\d^{12}(d_i)}}{\sum_{s=1}^k
    \frac{t_s}{\chi_\d^{8}(d_s)}} - 4 (1 + o_{\d,\Theta}(1))  \frac{ \sum_{i,s=1}^k t_i
    t_s \frac{\n (\chi_\d^2 (d_i)) \cdot \n (\chi_\d^2 (d_s))}{\chi_\d^{10}(d_i)
    \chi_\d^{10}(d_s)}}{\left[ \sum_{r=1}^k t_r \frac{1}{\chi_\d^8(d_r)}
     \right]^2} \nonumber \\ & + & o_{\d,\Theta}(1) \frac{ \sum_{i,s=1}^k t_i
    t_s \frac{|\n (\chi_\d^2 (d_i))| \; |\n (\chi_\d^2 (d_s))|}{\chi_\d^{10}(d_i)
    \chi_\d^{10}(d_s)}}{\left[ \sum_{r=1}^k t_r \frac{1}{\chi_\d^8(d_r)}
     \right]^2} \qquad \qquad \hbox{ in } M_{\s,\Theta}.
\end{eqnarray}
To have a further simplification of the last expression, it is
convenient to get rid of the cutoff functions $\chi_\d$ and
$\tilde{\chi}_\d$. In order to do this, we divide the set of
points $\{x_1, \dots, x_k\}$ in a suitable way. Since the number
$k$ is fixed, there exists $\hat{\d}$ and sets $\mathcal{C}_1,
\dots, \mathcal{C}_j$, $j \leq k$ with the following properties
\begin{equation}\label{eq:hatdCa}
    \left\{%
\begin{array}{ll}
    C_k^{-1} \d \leq \hat{\d} \leq \frac{\d}{16}; &  \\
    \mathcal{C}_1 \cup \dots \cup \mathcal{C}_j = \{x_1, \dots, x_k\}; &
    \\
    dist(x_i, x_s) \leq \hat{\d} & \hbox{ if } x_i, x_s \in \mathcal{C}_a; \\
    dist(x_i, x_s) > 4 \hat{\d} & \hbox{ if } x_i \in \mathcal{C}_a,
    x_s \in \mathcal{C}_b, a \neq b, \\
\end{array}%
\right.
\end{equation}
where $C_k$ is a positive constant depending only on $k$. Now we
define
$$
  \hat{\mathcal{C}}_a = \left\{ y \in M \; : \; dist(y, \mathcal{C}_a)
  \leq 2 \hat{\d} \right\}; \qquad T_a = \sum_{x_i \in \mathcal{C}_a}
  t_i, \qquad \qquad \hbox{ for } a = 1, \dots, j.
$$
By the definition of $\hat{\d}$ it follows that
\begin{equation}\label{eq:yinhC}
  \chi_\d(d_i(y)) = \tilde{\chi}_\d(d_i(y)) = d_i(y), \qquad \hbox{
  for } x_i \in \mathcal{C}_a \hbox{ and } y \in \hat{\mathcal{C}}_a,
\end{equation}
and
\begin{equation}\label{eq:youthC}
    \chi_\d(d_i(y)) \geq 2 \hat{\d}, \qquad \hbox{
  for } x_i \in \mathcal{C}_a \hbox{ and } y \not\in
  \hat{\mathcal{C}}_a.
\end{equation}
Furthermore one has
\begin{equation}\label{eq:disjhatC}
  \hat{\mathcal{C}}_a \cap \hat{\mathcal{C}}_b = \emptyset \hbox{
  for } a \neq b.
\end{equation}
From \eqref{eq:bettest} and \eqref{eq:youthC} it follows that
\begin{equation}\label{eq:estlaploutC}
  |\D \var_{\l,\s}(y)| \leq C_{\hat{\d}} \qquad \qquad \hbox{ for }
  y \in M \setminus \cup_{a=1}^j \hat{\mathcal{C}}_a.
\end{equation}
Therefore, by \eqref{eq:disjhatC}, it is sufficient to estimate
$\D \var_{\l,\s}$ inside each set $\hat{\mathcal{C}}_a$, where
\eqref{eq:yinhC} holds. We obtain immediately the following two
estimates, regarding the first terms in \eqref{eq:bettest}
\begin{equation}\label{eq:Ca1}
    \sum_{i=1}^k t_i \frac{\tilde{\chi}_\d^2(d_i)}{\chi_\d^{12}(d_i)} =
  \sum_{x_i \in \hat{\mathcal{C}}_a} \frac{t_i}{d_i^{10}} + O((1 - T_a)
  \hat{\d}^{-10}); \qquad \sum_{s=1}^k \frac{t_s}{\chi_\d^{8}(d_s)} =
  \sum_{x_s \in \hat{\mathcal{C}}_a} \frac{t_s}{d_s^{8}} + \ov{O}((1 - T_a)
  \hat{\d}^{-8}) \quad \hbox{ in } \hat{\mathcal{C}}_a.
\end{equation}
Here we have used the symbol $\ov{O}$ to denote a quantity such
that
$$
  \ov{O}(t) \geq C^{-1} t,
$$
where $C$ is large but fixed positive constant (which depends on
$k$, $M$, but not on $\d, \hat{\d}$, $\l$ and $(x_i)_i$). The same
dependence on the constants is understood when we use the symbol $O$
when it has as argument $(1-T_a)$, or its powers.

To estimate the second and the third term in the right-hand side of
\eqref{eq:bettest}, we use geodesic coordinates centered at some
point $y_a \in \hat{\mathcal{C}}_a$. With an abuse of notation, we
identify the points in $\mathcal{C}_a$ with their pre-image under
the exponential map. Using these coordinates, we find
$$
  \n (d_i(y))^2 = 2 (y - x_i) + o_\d(1) |y - x_i|, \qquad \hbox{ for }
  y \in \hat{\mathcal{C}}_a, \hbox{ and for } x_i \in \mathcal{C}_a,
$$
which implies
$$
    \frac{\n (\chi_\d^2 (d_i)) \cdot \n (\chi_\d^2 (d_s))}{\chi_\d^{10}(d_i)
    \chi_\d^{10}(d_s)} = 4 \frac{(y - x_i) \cdot (y - x_s)}{d_i^{10}
    d_s^{10}} + o_\d(1) \frac{1}{d_i^9 d_s^9}; \qquad \hbox{ for }
    y \in \hat{\mathcal{C}}_a \hbox{ and for } x_i, x_s \in \mathcal{C}_a.
$$
In particular, for $y \in \hat{\mathcal{C}}_a$, we get
\begin{eqnarray}\label{eq:Ca2}
    \sum_{i,s=1}^k t_i t_s \frac{\n (\chi_\d^2 (d_i)) \cdot \n
   (\chi_\d^2 (d_s))}{\chi_\d^{10}(d_i) \chi_\d^{10}(d_s)} & = &
   4 \sum_{x_i, x_s \in \mathcal{C}_a} t_i t_s \frac{(y - x_i) \cdot
   (y - x_s)}{d_i^{10} d_s^{10}} + o_\d(1) \sum_{x_i, x_s \in \mathcal{C}_a}
   \frac{t_i t_s}{d_i^9 d_s^9} \nonumber \\ & + & O((1 - T_a) \hat{\d}^{-9})
   \sum_{x_i \in \mathcal{C}_a} \frac{t_i}{d_i^9} + O((1 - T_a)^2
   \hat{\d}^{-18}).
\end{eqnarray}
We have also (still for $y \in \hat{\mathcal{C}}_a$)
\begin{eqnarray}\label{eq:Ca3}
    \sum_{i,s=1}^k t_i t_s \frac{|\n (\chi_\d^2 (d_i))| \; |\n
   (\chi_\d^2 (d_s))|}{\chi_\d^{10}(d_i) \chi_\d^{10}(d_s)} \leq
   4 \sum_{x_i, x_s \in \mathcal{C}_a} t_i t_s \frac{(y - x_i) \cdot
   (y - x_s)}{d_i^{10} d_s^{10}} + o_\d(1) \sum_{x_i, x_s \in \mathcal{C}_a}
   \frac{t_i t_s}{d_i^9 d_s^9}.
\end{eqnarray}
Hence from \eqref{eq:bettest}, \eqref{eq:Ca1}, \eqref{eq:Ca2} and
\eqref{eq:Ca3} we deduce (we are still working in the above
coordinates)
\begin{eqnarray*}
    \D \var_{\l,\s} (y) & = & 12 (1 + o_{\d,\Theta}(1)) \frac{\sum_{x_i \in
    \mathcal{C}_a} \frac{t_i}{d_i^{10}} + O((1 - T_a) \hat{\d}^{-10})}{\sum_{x_i
    \in \mathcal{C}_a} \frac{t_i}{d_i^{8}} + \ov{O}((1 - T_a) \hat{\d}^{-8})}
    \\ & - & 16 (1 + o_{\d,\Theta}(1)) \frac{ \left| \sum_{x_i \in \mathcal{C}_a}
    \frac{t_i (y-x_i)}{d_i^{10}} \right|^2 + o_\d(1) \left| \sum_{x_i \in
    \mathcal{C}_a} \frac{t_i}{d_i^9} \right|^2}{\left[ \sum_{x_i \in \mathcal{C}_a}
    \frac{t_i}{d_i^{8}} + \ov{O}((1 - T_a) \hat{\d}^{-8})\right]^2} \\
    & + & \frac{O((1 - T_a) \hat{\d}^{-9}) \sum_{x_i \in \mathcal{C}_a}
    \frac{t_i}{d_i^9} + O((1 - T_a)^2 \hat{\d}^{-18})}{\left[ \sum_{x_i
    \in \mathcal{C}_a} \frac{t_i}{d_i^{8}} + \ov{O}((1 - T_a) \hat{\d}^{-8})
    \right]^2}; \qquad \qquad y \in \hat{\mathcal{C}}_a.
\end{eqnarray*}
Using the inequality $ab \leq \e a^2 + C_\e b^2$ with $a = (1 - T_a)
\hat{\d}^{-9}$ and $b = \sum_{x_i \in \mathcal{C}_a}
    \frac{t_i}{d_i^9}$ we then find
\begin{eqnarray}\label{eq:LLCa}
    \D \var_{\l,\s} (y) & = & (1 + o_{\d,\Theta}(1)) \left[ 12
    \frac{\sum_{x_i \in \mathcal{C}_a} \frac{t_i}{d_i^{10}}}{\sum_{x_i
    \in \mathcal{C}_a} \frac{t_i}{d_i^{8}} + \ov{O}((1 - T_a)
    \hat{\d}^{-8})} - 16 \frac{ \left| \sum_{x_i \in \mathcal{C}_a}
    \frac{t_i (y-x_i)}{d_i^{10}} \right|^2}{\left[ \sum_{x_i \in
    \mathcal{C}_a} \frac{t_i}{d_i^{8}} + \ov{O}((1 - T_a) \hat{\d}^{-8})
    \right]^2} \right] \nonumber \\ & + & \frac{ (o_\d(1) + O(\e)) \left|
    \sum_{x_i \in \mathcal{C}_a} \frac{t_i}{d_i^9} \right|^2}{\left[ \sum_{x_i \in
    \mathcal{C}_a} \frac{t_i}{d_i^{8}} + \ov{O}((1 - T_a) \hat{\d}^{-8})
    \right]^2} + O(C_\e + 1) (1 - T_a)^2 \hat{\d}^{-2}; \qquad \qquad
    y \in \hat{\mathcal{C}}_a.
\end{eqnarray}
Now, given a large and fixed constant $\ov{C}$, we define the set
$\mathcal{B}^{\ov{C}}_a$ by
$$
  \mathcal{B}^{\ov{C}}_a = \left\{ y \in \hat{\mathcal{C}}_a \cap M_{\s,\Theta}
  \hbox{ s.t. if } x_i \in \mathcal{C}_a \hbox{ then } d_i(y) \leq \left( 1
  + \frac{1}{\ov{C}} \right) dist(y,\mathcal{C}_a) \hbox{ or } d_i(y) \geq
  \ov{C} dist(y,\mathcal{C}_a) \right\}.
$$
We start by characterizing the points belonging to the complement
of $\mathcal{B}^{\ov{C}}_a$ in $M_{\s,\Theta} \cap
\hat{\mathcal{C}}_a$. By definition, we have
\begin{equation}\label{eq:condcompl}
  y \in \left( M_{\s,\Theta} \cap \hat{\mathcal{C}}_a \right) \setminus
  \mathcal{B}^{\ov{C}}_a \quad \Rightarrow \quad \hbox{ there exists }
  x_i \in \mathcal{C}_a \hbox{ such that } d_i(y) \in \left( 1 +
  \frac{1}{\ov{C}}, \ov{C} \right) dist(y,\mathcal{C}_a).
\end{equation}
Given $y \in \left( M_{\s,\Theta} \cap \hat{\mathcal{C}}_a \right)
\setminus \mathcal{B}^{\ov{C}}_a$, we let $x_{\ov{i}}$ denote one of
its closest points in $\mathcal{C}_a$, and we let $x_{\ov{j}}$
denote one of the closest points in $\mathcal{C}_a$ to $y$, among
those which do not realize the infimum of the distance from $y$.
Then, since $dist(y,x_{\ov{i}}) < dist(y,x_{\ov{j}})$ and since
$dist(y,x_{\ov{j}}) < \ov{C} dist(y,x_{\ov{i}})$ (by
\eqref{eq:condcompl}), we clearly have
$$
  \frac{1}{\ov{C}} dist(y,x_{\ov{j}}) < dist(y,x_{\ov{i}}) <
  dist(y,x_{\ov{j}}),
$$
namely $y$ lies in an annulus centered at $x_{\ov{i}}$ whose radii
have a ratio equal to $\ov{C}$.

Now, fixing $x_{\ov{i}} \in \mathcal{C}_a$, we consider the
following set
$$
 \mathcal{D}_{\ov{i}} = \left\{ y \in \left( M_{\s,\Theta} \cap
 \hat{\mathcal{C}}_a \right) \setminus \mathcal{B}^{\ov{C}}_a
 \; : \; d_i(y) = dist(y, \mathcal{C}_a) \right\},
$$
namely the points $y$ in $\left( M_{\s,\Theta} \cap
\hat{\mathcal{C}}_a \right) \setminus \mathcal{B}^{\ov{C}}_a$ for
which $x_{\ov{i}}$ is the closest point to $y$ in $\mathcal{C}_a$.
Now, letting $y$ vary, there might be different points $x_{\ov{j}}$,
chosen as before, which do not realize the distance from $y$, but
anyway their number never exceeds $k$. This implies that
$\mathcal{D}_{\ov{i}}$ is contained in the union of at most $k$
annuli centered at $x_{\ov{i}}$ whose radii $c_{l,\ov{i}},
d_{l,\ov{i}}$ have uniformly bounded ratios, namely
\begin{equation}\label{eq:Di1}
    \mathcal{D}_{\ov{i}} \subseteq \cup_{l=1}^k \left( B_{d_{l,\ov{i}}}
  (x_{\ov{i}}) \setminus B_{c_{l,\ov{i}}} (x_{\ov{i}}) \right),
  \qquad \hbox{ with } d_{l,\ov{i}} \leq 2 \ov{C} c_{l,\ov{i}}.
\end{equation}
Clearly we also have
\begin{equation}\label{eq:Di2}
    \left( M_{\s,\Theta} \cap \hat{\mathcal{C}}_a \right) \setminus
  \mathcal{B}^{\ov{C}}_a = \cup_{x_{\ov{i}} \in \mathcal{C}_a}
  \mathcal{D}_{\ov{i}}.
\end{equation}
In $\mathcal{D}_{\ov{i}}$ there holds
$$
  \frac{t_i}{d_i^{10}} \leq \frac{1}{d_{\ov{i}}^2}
  \frac{t_i}{d_i^8}; \qquad \qquad \qquad \left| \sum_{x_i \in \mathcal{C}_a}
    \frac{t_i (y-x_i)}{d_i^{10}} \right| \leq \frac{1}{d_{\ov{i}}}
    \sum_{x_i \in \mathcal{C}_a} \frac{t_i}{d_i^{8}}.
$$
Then from \eqref{eq:LLCa} it follows that
\begin{equation}\label{eq:slfa}
    |\D \var_{\l,\s}| \leq C_{\d, \Theta, \e} \left( 1 + \frac{1}{d_{\ov{i}}^2}
  \right), \qquad \qquad \hbox{ in } \mathcal{D}_{\ov{i}}.
\end{equation}
Hence, from \eqref{eq:Di1} and \eqref{eq:Di2} using polar
coordinates we deduce that
\begin{eqnarray}\label{eq:D2Bac}
    \int_{\left( M_{\s,\Theta} \cap \hat{\mathcal{C}}_a \right)
    \setminus \mathcal{B}^{\ov{C}}_a} (\D \var_{\l,\s})^2 dV_g
    & \leq & \cup_{x_{\ov{i}} \in \mathcal{C}_a} \cup_{l=1}^k
    \int_{\left( B_{d_{l,\ov{i}}} (x_{\ov{i}}) \setminus
    B_{c_{l,\ov{i}}} (x_{\ov{i}}) \right)} C_{\d, \Theta, \e}
    \left( 1 + \frac{1}{d_{\ov{i}}^2} \right)^2 dV_g \nonumber \\ & \leq &
    C_{\d, \Theta, \e} \; k \; card(\mathcal{C}_a) \left(
    \log \frac{d_{l,\ov{i}}}{c_{l,\ov{i}}}
    + 1 \right) \leq C_{\d, \Theta, \e, \ov{C}}.
\end{eqnarray}
At this point, to estimate $\int_M (\D \var_{\l,\s})^2 dV_g$, it
only remains to consider the contribution inside
$\mathcal{B}^{\ov{C}}_a$.

In this set, we call $d_{a,min}$ the distance of $y$ from
$\mathcal{C}_a$, and $d_{a,out}$ the minimal distance of $y$ from
the points $x_i$ in $\mathcal{C}_a$ satisfying $d_i(y) \geq \ov{C}
dist(y, \mathcal{C}_a)$ (see the definition of
$\mathcal{B}^{\ov{C}}_a$). Therefore, setting
$$
  T_{a,in} = \sum_{x_i \in \mathcal{C}_a \; : \; d_i(y) \leq \left(
  1 + \frac{1}{\ov{C}} \right) dist(y,\mathcal{C}_a)} t_i,
$$
from \eqref{eq:LLCa} we obtain the estimate
\begin{eqnarray*}
    \D \var_{\l,x} (y) & = & 12 (1 + o_{\d,\Theta,\e,\ov{C}}(1)) \frac{
    \frac{T_{a,in}}{d_{a,min}^{10}} + \frac{O(T_a -
    T_{a,in})}{d_{a,out}^{10}}}{\frac{T_{a,in}}{d_{a,min}^8} + \frac{\ov{O}(T_a -
    T_{a,in})}{d_{a,out}^8} + \ov{O}((1 - T_a)\hat{\d}^{-8})} \\ & - & 16
    (1 + o_{\d,\Theta,\e,\ov{C}}(1)) \frac{ \left| \frac{T_{a,in}}{d_{a,min}^9} +
    \frac{O(T_a - T_{a,in})}{d_{a,out}^9} \right|^2}{\left[ \frac{T_{a,in}}{d_{a,min}^8} +
    \frac{\ov{O}(T_a - T_{a,in})}{d_{a,out}^8} + \ov{O}((1 - T_a)\hat{\d}^{-8})
    \right]^2} + C_{\d,\Theta,\e,\ov{C}}.
\end{eqnarray*}
Now we notice that for $y \in \mathcal{B}^{\ov{C}}_a$ the following
inequalities hold
$$
  \frac{T_{a,in}}{d_{a,min}^9} + \frac{O(T_a - T_{a,in})}{d_{a,out}^9} \leq \left(
  C_{\d, \Theta, \e, \ov{C}} + \frac{(1 + o_{\d,\Theta,\e,\ov{C}}(1))}{d_{a,min}} \right) \left(
  \frac{T_{a,in}}{d_{a,min}^8} +
  \frac{\ov{O}(T_a - T_{a,in})}{d_{a,out}^8} \right);
$$
$$
  \frac{T_{a,in}}{d_{a,min}^{10}} + \frac{O(T_a - T_{a,in})}{d_{a,out}^{10}} \leq \left(
  C_{\d, \Theta, \e, \ov{C}} + \frac{(1 + o_{\d,\Theta,\e,\ov{C}}(1))}{d_{a,min}^2} \right) \left(
  \frac{T_{a,in}}{d_{a,min}^8} +
  \frac{\ov{O}(T_a - T_{a,in})}{d_{a,out}^8} \right);
$$
$$
  \frac{T_{a,in}}{d_{a,min}^{10}} + \frac{O(T_a - T_{a,in})}{d_{a,out}^{10}}
  \geq \left( - C_{\d, \Theta, \e, \ov{C}} + \frac{(1 - o_{\d,\Theta,\e,\ov{C}}(1))}{d_{a,min}^2}
  \right) \left( \frac{T_{a,in}}{d_{a,min}^8} +
  \frac{\ov{O}(T_a - T_{a,in})}{d_{a,out}^8} \right).
$$
From the last four formulas and some elementary computations one
can deduce that
\begin{equation}\label{eq:elda}
    |\D \var_{\l,\s}| \leq C_{\d,\Theta,\e,\ov{C}} + 4 (1 + o_{\d,\Theta,\e,\ov{C}}(1))
  \frac{1}{d_{a,min}^2} \qquad \qquad \hbox{ in } \mathcal{B}_a^{\ov{C}}.
\end{equation}
We notice that, trivially
$$
  \mathcal{B}_a^{\ov{C}} = \cup_{x_i \in \mathcal{C}_a} \left(
  \mathcal{B}_a^{\ov{C}} \cap \{ y \; : \; d_i(y) = d_{a,min} \}
  \right).
$$
Therefore, recalling that $\frac{\Theta}{\l} \leq d_i(y) \leq
\hat{\d}$ for every $y \in \mathcal{B}_a^{\ov{C}}$, from the last
two formulas it follows that (the volume of the three-sphere is $2
\pi^2$)
\begin{eqnarray}\label{eq:L2ult}
    \int_{\mathcal{B}_a^{\ov{C}}} (\D \var_{\l,\s})^2 dV_g & \leq
    & \sum_{x_i \in \mathcal{C}_a} \int_{\mathcal{B}_a^{\ov{C}}
    \cap \{ y \; : \; d_i(y) = d_{a,min} \}} \left( C_{\d,\Theta,\e,\ov{C}}
    + 4 (1 + o_{\d,\Theta,\e,\ov{C}}(1)) \frac{1}{d_{a,min}^2} \right)^2
    dV_g \nonumber \\ & \leq & \sum_{x_i \in \mathcal{C}_a}
    \int_{B_{\hat{\d}}(x_i) \setminus B_{\frac{\Theta}{\l}}(x_i)}
    \left( C_{\d,\Theta,\e,\ov{C}} + 4 (1 + o_{\d,\Theta,\e,\ov{C}}(1))
  \frac{1}{d_{a,min}^2} \right)^2 dV_g \\ & \leq & \nonumber
    card(\mathcal{C}_a) \left( 32 \pi^2 (1 + o_{\d,\Theta,\e,\ov{C}}(1))
    \log \frac{\hat{\d} \l}{\Theta} + C_{\d,\Theta,\e,\ov{C}} \right)
    \\ & \leq & card(\mathcal{C}_a) 32 \pi^2 (1 + o_{\d,\Theta,\e,\ov{C}}(1))
    \log \l + C_{\d,\Theta,\e,\ov{C}} \nonumber .
\end{eqnarray}
From \eqref{eq:intL2Bi}, \eqref{eq:estlaploutC}, \eqref{eq:D2Bac}
and \eqref{eq:L2ult}, considering all the sets
$\hat{\mathcal{C}}_a$ and the complement of their union, we
finally deduce
\begin{equation}\label{eq:La2m}
    \int_M (\D \var_{\l,\s})^2 dV_g \leq 32 \pi^2 k (1 +
    o_{\d,\Theta,\e,\ov{C}}(1)) \log \l +
    C_{\hat{\d},\e,\ov{C},\Theta}.
\end{equation}
Fixing the values of $\ov{C}, \Theta$ (large, depending on $\d$) and
of $\e$ (small, depending on $\d$), we obtain the estimate of the
term involving the squared Laplacian.

Next, we estimate the term $\int_M |\n \var_{\l,\s}|^2 dV_g$. It
could be possible to proceed using $L^p$ estimates on $\var_{\l,\s}
- \ov{\var_{\l,\s}}$ and interpolation, but having the computations
for the Laplacian at hand, it is convenient to work directly. From
\eqref{eq:njvar}, one finds first the following pointwise estimate
$$
  |\n \var_{\l,\s}| \leq \frac{C}{\l},
$$
which implies, similarly as before
\begin{equation}\label{eq:intL2Bi2}
    \int_{\cup_{i=1}^k B_{\frac{\Theta}{\l}}(x_i)} |\n \var_{\l,\s}|^2
    dV_g \leq C \frac{\Theta^4}{\l^2}.
\end{equation}
On the other hand, in the set $M_{\s,\Theta}$, using
\eqref{eq:appth} and reasoning as above we obtain
$$
    \n \var_{\l,\s} = - (1 + o_{\d,\Theta}(1))
    \frac{ \sum_i t_i \frac{\n (\chi_\d^2(d_i))}{\chi_\d^{10}(d_i)}}{\sum_s
\frac{t_s}{\chi_\d^{8}(d_s)}} + o_{\d,\Theta}(1) \frac{ \sum_i t_i
\frac{|\n (\chi_\d^2(d_i))|}{\chi_\d^{10}(d_i)}}{\sum_s
\frac{t_s}{\chi_\d^{8}(d_s)}}.
$$
Taking the square we get
\begin{eqnarray}\label{eq:bettest2}
    |\n \var_{\l,\s}|^2 & \leq  & (1 + o_{\d,\Theta}(1))
    \frac{ \sum_{i,s} t_i t_s \frac{\n (\chi_\d^2(d_i)) \cdot \n
    (\chi_\d^2(d_s))}{\chi_\d^{10}(d_i) \chi_\d^{10}(d_s)}}{\left[\sum_s
\frac{t_s}{\chi_\d^{8}(d_s)}\right]^2} + o_{\d,\Theta}(1) \left(
\frac{ \sum_i t_i \frac{|\n
(\chi_\d^2(d_i))|}{\chi_\d^{10}(d_i)}}{\sum_s
\frac{t_s}{\chi_\d^{8}(d_s)}} \right)^2.
\end{eqnarray}
Using \eqref{eq:Ca2} and \eqref{eq:Ca3} we deduce (working as
before in geodesic coordinates)
\begin{eqnarray*}
    |\n \var_{\l,\s}|^2 (y) & = & 4 (1 + o_{\d,\Theta,\e}(1)) \frac{ \left|
    \sum_{x_i \in \mathcal{C}_a}
    \frac{t_i (y-x_i)}{d_i^{10}} \right|^2 + o_{\d,\Theta,\e}(1)
    \left| \sum_{x_i \in \mathcal{C}_a} \frac{t_i}{d_i^9} \right|^2}{\left[
    \sum_{x_i \in \mathcal{C}_a} \frac{t_i}{d_i^{8}} + \ov{O}((1 - T_a)
    \hat{\d}^{-8})\right]^2} \\ & + & \frac{C_{\d,\Theta,\e} O((1 - T_a)^2
    \hat{\d}^{-18})}{\left[ \sum_{x_i \in \mathcal{C}_a} \frac{t_i}{d_i^{8}}
    + \ov{O}((1 - T_a) \hat{\d}^{-8})\right]^2}, \qquad \qquad y \in
   \hat{\mathcal{C}}_a.
\end{eqnarray*}
Reasoning as for \eqref{eq:slfa} and \eqref{eq:elda}, one then
finds
$$
  |\n \var_{\l,\s}|^2 \leq C_{\d,\Theta,\e,\ov{C}} \left( 1 +
  \frac{1}{d_{a,min}^2} \right) \qquad \qquad \hbox{ in }
  \hat{\mathcal{C}}_a \cap M_{\s,\Theta},
$$
which implies
$$
  \int_{\hat{\mathcal{C}}_a} |\n \var_{\l,\s}|^2 dV_g \leq
  C_{\d,\Theta,\e,\ov{C}}.
$$
On the other hand, we have also
$$
  |\n \var_{\l,\s}(y)|^2 \leq C_{\hat{\d}} \qquad \qquad \hbox{ for }
  y \in M \setminus \cup_{a=1}^j \hat{\mathcal{C}}_a.
$$
Therefore from the last two formulas we deduce
\begin{equation}\label{eq:estfingrad2}
  \int_{M} |\n \var_{\l,\s}|^2 dV_g \leq
  \hat{C}_{\d,\Theta,\e,\ov{C}}.
\end{equation}
From \eqref{eq:Pguv} it follows that
$$
  \langle P_g \var_{\l,\s}, \var_{\l,\s} \rangle \leq \int_M (\D
  \var_{\l,\s})^2 dV_g + C \int_M |\n \var_{\l,\s}|^2 dV_g
$$
Hence, from \eqref{eq:La2m} and \eqref{eq:estfingrad2} we finally
obtain, fixing as before the values of the constants $\Theta,\e$ and
$\ov{C}$
$$
  \langle P_g \var_{\l,\s}, \var_{\l,\s} \rangle \leq 32 k \pi^2
  (1 + o_\d(1)) \log \l + C_\d.
$$
This concludes the proof.
\end{pfnb}

\

\begin{pfnb} {\bf of Lemma \ref{l:arho}} \quad
If $P_g$ is non-negative, for $8 (k+1) \pi^2 > \rho' \geq \rho > 8 k
\pi^2$ (resp. for $8 \pi^2 > \rho' \geq \rho$) we clearly have
$$
  \frac{II_\rho(u)}{\rho} - \frac{II_{\rho'}(u)}{\rho'} = \left(
  \frac 1\rho - \frac{1}{\rho'} \right) \langle P_g u, u \rangle
  \geq 0,
$$
and the conclusion follow immediately taking $C = 0$. Therefore from
now on we consider the case in which $P_g$ possesses some negative
eigenvalues. The last formula in this case yields
\begin{equation}\label{eq:ineqP-}
    \frac{II_{\rho'}(u)}{\rho'} \leq \frac{II_\rho(u)}{\rho} -
    \frac{(\rho' - \rho)}{\rho \rho'} \langle P_g \hat{u}, \hat{u}
    \rangle,
\end{equation}
where $\hat{u}$ is the $V$-component of $u$, see \eqref{eq:hatv1k}.

Fixing $\rho \in [1 - \rho_0, 1 + \rho_0]$ and $\e > 0$, we consider
a map $\pi_{\rho,\e} \in \ov{\Pi}_{\ov{S},\ov{\l}}$ such that
\begin{equation}\label{eq:suppre}
    \sup_{m \in \widehat{A_{k,\ov{k}}}} II_{\rho}
    (\pi_{\rho,\e}(m)) < \ov{\Pi}_\rho + \e.
\end{equation}
We can also assume that each element of the form $u =
\pi_{\rho,\e}(m)$ satisfies the normalization condition $\int_M e^{4
u} dV_g = 1$. Now, considering the $V$-part $\hat{u}$ of all these
elements, we fix three numbers $\th > 0$ (small, depending on
$\pi_{\rho,\e}$), and $C_0, C_1 > 0$ (depending on $M$ and $\rho_0$,
with $C_1 \gg C_0 \gg 1$), and we define a new map
$\tilde{\pi}_{\rho,\e}$ in the following way
\begin{equation}\label{eq:tpier}
  \tilde{\pi}_{\rho,\e}(m) = \pi_{\rho,\e}(m) + \eta_\th(m)
  \tilde{\eta}(\widehat{\pi_{\rho,\e}(m)})
  \widehat{\pi_{\rho,\e}(m)}; \qquad \quad m \in
  \widehat{A_{k,\ov{k}}}.
\end{equation}
Here $\widehat{\pi_{\rho,\e}(m)}$ denotes the $V$-component of
$\pi_{\rho,\e}(m)$ (see Section \ref{s:pr}), the function
$\eta_\th(m)$, $m = (m_1, t) \in A_{k,\ov{k}} \times [0,1]$, is
defined as
$$
  \eta_\th(m) = \left\{%
\begin{array}{ll}
    1, & \hbox{ for } t \in [0, 1- \th]; \\
    \frac{1}{\th} (1-t), & \hbox{ for } t \in [1-\th,1], \\
\end{array}%
\right.
$$
and $\tilde{\eta}(\widehat{\pi_{\rho,\e}(m)})$ is given by
$$
  \tilde{\eta}(\widehat{\pi_{\rho,\e}(m)}) = \left\{%
\begin{array}{ll}
    0, & \hbox{ for } \|\widehat{\pi_{\rho,\e}(m)}\| \in [0,C_0]; \\
    \frac{1}{C_1 - C_0} (\|\widehat{\pi_{\rho,\e}(m)}\| - C_0), &
    \hbox{ for } \| \widehat{\pi_{\rho,\e}(m)} \| \in [C_0, C_1]; \\
    1, & \hbox{ for } \| \widehat{\pi_{\rho,\e}(m)} \| \geq C_1. \\
\end{array}%
\right.
$$

When $\eta_\th(m) = 1$, by the normalization of $\pi_{\rho,\e}$ we
have the following upper bound on
$II_{\rho}(\tilde{\pi}_{\rho,\e}(m))$
\begin{eqnarray}\label{eq:prima} \nonumber
  II_{\rho}(\tilde{\pi}_{\rho,\e}) & = & \langle P_g
  \pi_{\rho,\e}, \pi_{\rho,\e} \rangle + \left(2 \tilde{\eta}(\widehat{\pi_{\rho,\e}})
  + (\tilde{\eta}(\widehat{\pi_{\rho,\e}}))^2 \right) \langle P_g
  \widehat{\pi_{\rho,\e}}, \widehat{\pi_{\rho,\e}} \rangle \\ \nonumber \\ & +
  & 4 \rho \int_M Q_g ( \pi_{\rho,\e} +
  \tilde{\eta}(\widehat{\pi_{\rho,\e}})
  \widehat{\pi_{\rho,\e}} ) d V_g - 4 \rho k_P \int_M e^{4 \pi_{\rho,\e} +
  4 \tilde{\eta}(\widehat{\pi_{\rho,\e}}) \widehat{\pi_{\rho,\e}}} d V_g \\ \nonumber \\ & \leq &
  II_{\rho}(\pi_{\rho,\e}) + \left(2 \tilde{\eta}(\widehat{\pi_{\rho,\e}})
  + (\tilde{\eta}(\widehat{\pi_{\rho,\e}}))^2 \right) \langle P_g
  \widehat{\pi_{\rho,\e}}, \widehat{\pi_{\rho,\e}} \rangle +
  \tilde{C}_0 \tilde{\eta}(\widehat{\pi_{\rho,\e}})
  \left\| \widehat{\pi_{\rho,\e}} \right\|, \nonumber
\end{eqnarray}
where $\tilde{C}_0$ is a fixed constant depending only on $M$ and on
$\rho_0$.

Since $\widehat{\pi_{\rho,\e}}$ belongs to the space $V$, where
$P_g$ is negative-definite, if $C_0$ is sufficiently large
(depending only $\tilde{C}_0$ which, in turn, depends only on $M$
and $\rho_0$), then one has
\begin{equation}\label{eq:minor}
    \left(2 \tilde{\eta}(\widehat{\pi_{\rho,\e}})
  + (\tilde{\eta}(\widehat{\pi_{\rho,\e}}))^2 \right) \langle P_g
  \widehat{\pi_{\rho,\e}}, \widehat{\pi_{\rho,\e}} \rangle +
  \tilde{C}_0 \tilde{\eta}(\widehat{\pi_{\rho,\e}})
  \left\| \widehat{\pi_{\rho,\e}} \right\| \leq 0 \qquad \hbox{
  for } \left\| \widehat{\pi_{\rho,\e}}(m) \right\| \geq C_0.
\end{equation}
Having fixed this value of $C_0$, from \eqref{eq:ineqP-} and the
fact that $\tilde{\eta}(\widehat{\pi_{\rho,\e}(m)}) = 0$ for
$\|\tilde{\eta}(\widehat{\pi_{\rho,\e}(m)})\| \leq C_0$ it follows
that
\begin{equation}\label{eq:c0c0c0}
    \frac{II_{\rho'}(\tilde{\pi}_{\rho,\e})}{\rho'} \leq
    \frac{II_{\rho}(\pi_{\rho,\e})}{\rho} - \frac{\rho' - \rho}{\rho
    \rho'} \langle P_g \widehat{{\pi}_{\rho,\e}},
    \widehat{{\pi}_{\rho,\e}} \rangle \leq \frac{\ov{\Pi}_\rho +
    \e}{\rho} + \hat{C}_0 (\rho' - \rho); \qquad
    \|\widehat{{\pi}_{\rho,\e}}(m)\| \leq C_0, \eta_\th(m) = 1,
\end{equation}
where $\hat{C}_0$ depends only on $M$ and $\rho_0$.

Now we fix also the value of $C_1$. We choose $\rho_0$ so small and
$C_1 > 0$ (depending only on $M$ and $\rho_0$) so large that
\begin{equation}\label{eq:c1c1}
  \frac{3}{\rho} \left( 1 - \frac 43 \frac{\rho'-\rho}{\rho'}
  \right) \langle P_g \hat{v}, \hat{v} \rangle + \frac{\tilde{C}_0}{\rho}
  \|\hat{v}\| \leq \langle P_g \hat{v}, \hat{v} \rangle \leq - 2
  \ov{L} - L \qquad \hbox{ for all } \hat{v} \in V \hbox{ with }
  \|\hat{v}\| \geq C_1,
\end{equation}
where $L$ and $\ov{L}$ are the constants given in
\eqref{eq:min-maxrho} and \eqref{eq:ovlovl}. From \eqref{eq:ineqP-},
\eqref{eq:prima} and \eqref{eq:minor} we immediately find (still for
$\eta_\th(m) = 1$)
\begin{equation}\label{eq:c1c1c1}
    \frac{II_{\rho'}(\tilde{\pi}_{\rho,\e})}{\rho'} \leq
    \frac{II_{\rho}(\pi_{\rho,\e})}{\rho} - \frac{\rho' - \rho}{\rho
    \rho'} \langle P_g \widehat{{\pi_{\rho,\e}}},
    \widehat{{\pi_{\rho,\e}}} \rangle \leq \frac{\ov{\Pi}_\rho +
    \e}{\rho} + \hat{C}_1 (\rho' - \rho); \qquad
    C_0 \leq \|\widehat{{\pi_{\rho,\e}(m)}}\| \leq C_1,
\end{equation}
where $\hat{C}_1$ depends only on $M$ and $\rho_0$.

By \eqref{eq:ineqP-} and \eqref{eq:prima}, since
$\tilde{\eta}(\widehat{\pi_{\rho,\e}}) = 1$ when
$\|\widehat{\pi_{\rho,\e}}\| \geq C_1$ (which implies
$\widehat{\tilde{\pi}} = 2 \hat{\pi}$), we obtain
\begin{equation}\label{eq:iirhotp}
    \frac{II_{\rho}(\tilde{\pi}_{\rho,\e})}{\rho'} \leq
  \frac{II_{\rho}(\pi_{\rho,\e})}{\rho} + \frac{3}{\rho}
  \left( 1 - \frac 43 \frac{\rho'-\rho}{\rho'} \right) \langle P_g
  \widehat{\pi_{\rho,\e}}, \widehat{\pi_{\rho,\e}} \rangle +
  \frac{\tilde{C}_0}{\rho} \left\| \widehat{\pi_{\rho,\e}}
  \right\|; \qquad \|\widehat{\pi_{\rho,\e}(m)}\| \geq C_1,
  \eta_\th(m) = 1.
\end{equation}
Then \eqref{eq:c1c1} implies (see \eqref{eq:ovlovl} and
\eqref{eq:suppre})
\begin{equation}\label{eq:c2c2}
    \frac{II_{\rho'}(\tilde{\pi}_{\rho,\e})}{\rho'} \leq
  \frac{\ov{\Pi}_\rho}{\rho}, \qquad \qquad \hbox{ for }
  \|\hat{\pi}\| \geq C_1.
\end{equation}
From \eqref{eq:c0c0c0}, \eqref{eq:c1c1c1} and \eqref{eq:c2c2} we
deduce
\begin{equation}\label{eq:c3c3}
  \frac{II_{\rho'}(\tilde{\pi}_{\rho,\e})}{\rho'} \leq
  \frac{\ov{\Pi}_\rho + \e}{\rho} + (\hat{C}_0 + \hat{C}_1) (\rho' -
  \rho) \quad \qquad \hbox{ for } \eta_\th(m) = 1.
\end{equation}
Therefore it remains to consider the case in which $\eta_\th(m) \neq
1$, namely for $t > 1 - \th$ (recall that $m = (m_1,t)$ with $m_1
\in A_{k,\ov{k}}$). This is where the choice of $\th$ enters.
Reasoning as for \eqref{eq:prima} we find
\begin{eqnarray}\label{eq:prima2} \nonumber
  II_{\rho'}(\tilde{\pi}_{\rho,\e}) \leq
  II_{\rho'}(\pi_{\rho,\e}) + 2 \eta_\th(m) \tilde{\eta}(\widehat{\pi_{\rho,\e}})
  \langle P_g \widehat{\pi_{\rho,\e}}, \widehat{\pi_{\rho,\e}} \rangle +
  \tilde{C}_0 \eta_\th(m) \tilde{\eta}(\widehat{\pi_{\rho,\e}})
  \left\| \widehat{\pi_{\rho,\e}} \right\|. \nonumber
\end{eqnarray}
Recall that the map $\pi_{\rho,\e}$ belongs to
$\Pi_{\ov{S},\ov{\l}}$, and hence it satisfies
$\pi_{\rho,\e}(t,\cdot) \to \Phi_{\ov{S},\ov{\l}}(\cdot)$ in
$C^0(A_{k,\ov{k}})$. Since $\rho$ is varying in the small interval
$[ 1 - \rho_0, 1 + \rho_0]$, we have estimates of the form
\eqref{eq:IItest} (with $\rho k_P$ replacing $k_P$) uniformly for
$\rho$ in this interval. Thus from the last formula we deduce that,
for $\eta_\th(m) < 1$
\begin{eqnarray*}
  II_{\rho'}(\tilde{\pi}_{\rho,\e}(m)) & \leq &
  II_{\rho'}(\Phi_{\ov{S},\ov{\l}}(m_1)) + o_\th(1) +
  2 \eta_\th(m) \tilde{\eta}(\widehat{\pi_{\rho,\e}})
  \langle P_g \widehat{\pi_{\rho,\e}}, \widehat{\pi_{\rho,\e}} \rangle +
  \tilde{C}_0 \eta_\th(m) \tilde{\eta}(\widehat{\pi_{\rho,\e}})
  \left\| \widehat{\pi_{\rho,\e}} \right\| \\ \\ & \leq & \left( 32 k \pi^2
  - 4 \rho' k_P + o_\d(1) \right) \log \l - |\l_{\ov{k}}| |s|^2 \ov{S}^2
  + O(|s| \ov{S}) + C_\d + O(1) + o_\th(1) \\ \\ & + & 2 \eta_\th(m)
  \tilde{\eta}(\widehat{\pi_{\rho,\e}}) \langle P_g \widehat{\pi_{\rho,\e}},
  \widehat{\pi_{\rho,\e}} \rangle + \tilde{C}_0 \eta_\th(m)
  \tilde{\eta}(\widehat{\pi_{\rho,\e}}) \left\| \widehat{\pi_{\rho,\e}}
  \right\| \\ \\ & \leq & \left( 32 k \pi^2
  - 4 \rho' k_P + o_\d(1) \right) \log \l - |\l_{\ov{k}}| |s|^2 \ov{S}^2
  + O(|s| \ov{S}) + C_\d + O(1) < - \frac 32 L,
\end{eqnarray*}
if $L$ is chosen sufficiently large (see \eqref{eq:min-maxrho}) and
$\th$ is chosen sufficiently small. Now the conclusion follows from
\eqref{eq:c3c3} and the last estimate.
\end{pfnb}

\end{document}